\documentclass[reqno, 11pt, a4paper]{amsart}

\usepackage[T1]{fontenc}
\usepackage{amssymb}
\usepackage{hyperref}
\usepackage{graphicx}
\usepackage{tikz-cd}
\usepackage{bbm}

\begin{document}

\title{On the tautological rings of $\mathcal{M}_{g, 1}$ and its universal Jacobian}
\author{Qizheng Yin}
\dedicatory{\rm with an appendix by {\sc Li Ma}}
\address{Université Pierre et Marie Curie Paris VI, 4 place Jussieu, 75252 Paris cedex 05, France \newline \indent University of Amsterdam, Korteweg-de Vries Institute for Mathematics, P.O. Box 94248, 1090 GE Amsterdam, The Netherlands}
\email{qyin@uva.nl}
\address{Université Pierre et Marie Curie Paris VI, 4 place Jussieu, 75252 Paris cedex 05, France}
\email{lima@math.jussieu.fr}
\date{\today}

\begin{abstract}
We give a new method of producing relations in the tautological ring $\mathcal{R}(\mathcal{M}_{g, 1})$, using the $\mathfrak{sl}_2$-action on the Chow ring of the universal Jacobian. With these relations, we prove that $\mathcal{R}(\mathcal{M}_{g, 1})$ is generated by $\kappa_1, \ldots, \kappa_{\lfloor g/3 \rfloor}$ and $\psi$. Our computation shows that Faber's conjectures for $\mathcal{M}_{g, 1}$ are true for $g \leq 19$. Further, by pushing relations forward to $\mathcal{M}_g$, we obtain a new proof of Faber's conjectures (for $\mathcal{M}_g$) for $g \leq 23$. For $g = 24$, our method recovers all the Faber-Zagier relations. We also give an algebraic proof of an identity of Morita.

\end{abstract}
\maketitle

\medskip
\section*{\bf Introduction}
\medskip

\subsection*{\it Moduli side} --- Denote by $\mathcal{M}_{g, 1}$ the moduli space of smooth $1$-pointed curves of genus $g \geq 2$. Mumford and Faber initiated the study of the `tautological ring', which is defined to be a $\mathbb{Q}$-subalgebra of the Chow ring of the moduli space (with $\mathbb{Q}$-coefficients) generated by certain geometrically constructed classes. In the context of $\mathcal{M}_{g, 1}$, the generators are the classes $\{\kappa_i\}_{i \geq 0}$ that contain information about the curves, and $\psi$ that contains information about the marked point. We denote this tautological ring by $\mathcal{R}(\mathcal{M}_{g, 1})$.

An analogy of Faber's conjectures ({\it cf.}~\cite{Fab99}) in this case predicts the following: that $\mathcal{R}(\mathcal{M}_{g, 1})$ is Gorenstein with socle in codimension $g - 1$, that it is generated by $\kappa_1, \ldots, \kappa_{\lfloor g/3 \rfloor}$ and $\psi$, and that there are no relations between these classes in codimension $\leq \lfloor g/3 \rfloor$ ({\it cf.}~Conjectures~\ref{faber}). The difficulty of solving these conjectures is to find sufficiently many relations between tautological classes.

\subsection*{\it Jacobian side} --- Denote by $\pi \colon J \to \mathcal{M}_{g, 1}$ the universal Jacobian over $\mathcal{M}_{g, 1}$. Using the marked point we can embed the universal curve $C$ in $J$. There is also a notion of `tautological ring' for $J$, which is defined to be the smallest $\mathbb{Q}$-subalgebra of the Chow ring of $J$ that contains the curve class $[C]$, and that is stable under certain operations ({\it cf.}~Definition~\ref{defTJ}). We denote this tautological ring by $\mathcal{T}(J)$. 

In several papers \cite{Pol05}, \cite{Pol07} and \cite{Pol07b}, Polishchuk studied the $\mathfrak{sl}_2$-action on the Chow ring of Jacobians. By applying his results to our setting, we obtain that $\mathcal{T}(J)$ is generated by an explicit finite set of classes $\{p_{i, j}\}$, together with $\psi$. We also show that the operator $f \in \mathfrak{sl}_2$ acts on $\mathcal{T}(J)$ via an explicit differential operator $\mathcal{D}$ on the polynomial ring in $\{p_{i, j}\}$ and $\psi$ ({\it cf.}~Theorem~\ref{structure}). Further, the ring $\mathcal{R}(\mathcal{M}_{g, 1})$ can be identified as a $\mathbb{Q}$-subalgebra of $\mathcal{T}(J)$ ({\it cf.}~Corollary~\ref{subring}).

Polishchuk's approach in \cite{Pol05} provides a powerful method to produce relations in $\mathcal{T}(J)$: take any polynomial in $\{p_{i, j}\}$ and $\psi$ that vanishes for `obvious' reasons ({\it cf.}~Section~\ref{relation}), then apply the operator $\mathcal{D}$ one or several times. The resulting polynomial should vanish as well. In this way we obtain a huge space of relations that are simply dictated by the $\mathfrak{sl}_2$-action. By restricting everything to the subalgebra $\mathcal{R}(\mathcal{M}_{g, 1})$, we obtain relations there, too.

\subsection*{\it Main results} --- Here we list the main results of this note.

\smallskip
\noindent\makebox[\leftmargin][r]{(i) }Using our relations, we prove that $\mathcal{R}(\mathcal{M}_{g, 1})$ is generated by the classes $\kappa_1, \ldots, \kappa_{\lfloor g/3 \rfloor}$ and $\psi$ ({\it cf.}~Theorem~\ref{generation}). 

\smallskip
\noindent\makebox[\leftmargin][r]{(ii) }With the help of a computer, we confirm that Faber's conjectures for $\mathcal{M}_{g, 1}$ are true for $g \leq 19$ ({\it cf.}~Theorem~\ref{g19}). From $g = 20$ on, our relations are not enough to conclude that $\mathcal{R}(\mathcal{M}_{g, 1})$ is Gorenstein ({\it cf.}~Table~\ref{tab}).

\smallskip
\noindent\makebox[\leftmargin][r]{(iii) }By pushing forward to $\mathcal{M}_g$, we reprove Faber's conjectures (for $\mathcal{M}_g$) for $g \leq 23$ ({\it cf.}~Theorem~\ref{g23}). For $g = 24$ and so on, computation gives the same set of relations as the Faber-Zagier relations ({\it cf.}~Table~\ref{tab}). 

\smallskip
\noindent\makebox[\leftmargin][r]{(iv) }We also give an algebraic proof of an identity obtained by Morita ({\it cf.}~Theorem~\ref{morita}).

\smallskip
We refer to Sections \ref{socle} and \ref{sequence} for some other results and insights. Beyond all these results, our approach has many advantages compared to previous ones. From a theoretical perspective, it gives a clean and uniform treatment of Faber's conjectures, which converts a geometric problem into a combinatorial problem. With the help of the `Dutch house' introduced in Section~\ref{dutchhouse}, many complicated facts become obvious. On the practical side, it produces huge quantities of relations, considerably more than other methods.

Finally, the nature of our approach (using the $\mathfrak{sl}_2$-action as source of relations) also suggests that these might be the only relations we can ever find ({\it cf.}~Conjecture~\ref{conj}).

\subsection*{\it Notation and conventions} --- (i) We fix a base field $k$ of arbitrary characteristic.

\smallskip
\noindent\makebox[\leftmargin][r]{(ii) }For a smooth variety $X$ over $k$, we denote by ${\rm CH}(X)$ the Chow ring of $X$ with $\mathbb{Q}$-coefficients. We write $\mathbbm{1}_X$ for the fundamental class of $X$, which is the unit for the intersection product of ${\rm CH}(X)$.

\smallskip
\noindent\makebox[\leftmargin][r]{(iii) }Let $S$ be a smooth connected variety over $k$, and let $p \colon C \to S$ be a smooth relative curve of genus $g \geq 2$ over $S$. Denote by $K$ the first Chern class of $\Omega^1_{C/S}$, which belongs to ${\rm CH}^1(C)$. Following Mumford, we define for $i \geq 0$ the classes
\begin{equation*}
\kappa_i := p_*(K^{i + 1}) \in {\rm CH}^i(S).
\end{equation*}
We have $\kappa_0 = (2g - 2) \cdot \mathbbm{1}_S$, and it is often convenient to extend this definition by adding $\kappa_{-1} = 0$.

\smallskip
\noindent\makebox[\leftmargin][r]{(iv) }Assume that $p \colon C \to S$ admits a section (marked point) $x_0 \colon S \to C$. We define the class
\begin{equation*}
\psi := x_0^*(K) \in {\rm CH}^1(S).
\end{equation*}
Note that if we write $[x_0] := \big[x_0(S)\big] \in {\rm CH}^1(C)$, we have $x_0^*\big([x_0]\big) = - \psi$ by adjunction. By abuse of notation, we shall write the same $\psi$ for all pull-backs of $\psi$ to varieties over $S$ ({\it e.g.}~the relative curve $C$ or the relative Jacobian $J$).

\smallskip
\noindent\makebox[\leftmargin][r]{(v) }Following Mumford and Faber, we consider the $\mathbb{Q}$-subalgebra of ${\rm CH}(S)$ generated by the geometrically constructed classes $\{\kappa_i\}_{i \geq 0}$ and $\psi$. We call it the {\it tautological ring\/} of $S$, and denote it by $\mathcal{R}(S)$.

\smallskip
\noindent\makebox[\leftmargin][r]{(vi) }For $g \geq 2$, let $\mathcal{M}_{g, 1}$ be the moduli stack of smooth $1$-pointed curves of genus $g$. It admits a finite cover by a smooth quasi-projective variety. Since we work with $\mathbb{Q}$-coefficients, the Chow ring ${\rm CH}(\mathcal{M}_{g, 1})$ can be easily defined via the cover. Similarly, we define the classes $\{\kappa_i\}$ and $\psi$ in ${\rm CH}(\mathcal{M}_{g, 1})$, and also the tautological ring $\mathcal{R}(\mathcal{M}_{g, 1})$. In principle, we may regard $\mathcal{M}_{g, 1}$ as a smooth connected variety when talking about Chow groups with $\mathbb{Q}$-coefficients.

\smallskip
\noindent\makebox[\leftmargin][r]{(vii) }We write $\mathfrak{sl}_2 = \mathbb{Q} \cdot e + \mathbb{Q} \cdot f + \mathbb{Q} \cdot h$, with $[e, f] = h$, $[h, e] = 2e$ and $[h, f] = -2f$.

\subsection*{\it Acknowledgements} --- I am deeply indebted to my thesis advisor Ben~Moonen, who introduced me to this subject and encouraged me throughout this project. Without him I would have given up long time ago. I am especially grateful to my colleague Li~Ma, whose brilliant computer program made this work much more meaningful. I thank Claire~Voisin (my other thesis advisor), Carel~Faber, Mehdi~Tavakol and Olof~Bergvall for useful discussions and feedbacks. The whole project is originally motivated by Robin~de~Jong's talk at the PCMI in summer 2011, and many ideas are due to Alexander~Polishchuk.

I wish to thank my parents and my fiancée Zhiyue~Zhou for their constant support. Finally I am grateful to a childhood friend who inspired me on New Year's Eve 2011.

\medskip
\section{\bf The Beauville decomposition, the Dutch house and the \texorpdfstring{$\mathfrak{sl}_2$}{sl2}-action}
\medskip

\noindent In this section we explain two basic ingredients for our study. One is the classical result of Beauville \cite{Bea86} on the decomposition of the Chow ring of abelian varieties, later generalized to abelian schemes by Deninger and Murre \cite{DM91}. Following Moonen \cite{Moo09},   we use an illustration called the `Dutch house' to explain this decomposition. The other ingredient is the $\mathfrak{sl}_2$-action on the Chow ring of abelian schemes induced by the Lefschetz operator, as studied by Künnemann \cite{Kün93}. For the latter we shall focus on the case of the relative Jacobian associated to a relative $1$-pointed curve. Polishchuk \cite{Pol07b} has shown that in this situation, the $\mathfrak{sl}_2$-action can be described geometrically.

\subsection{\it The Beauville decomposition} --- Take $S$ a smooth connected variety of dimension $d$ over $k$. Let $\pi \colon A \to S$ be an abelian scheme of relative dimension $g$, together with a principal polarization $\lambda \colon A \xrightarrow{\sim} A^t$. We use $\lambda$ to identify $A$ with its dual abelian scheme $A^t$. 

Denote by $[n] \colon A \to A$ the multiplication by $n$. According to Deninger and Murre \cite{DM91},  the relative Chow motive of $A$ over $S$ (denoted by $R(A/S)$) admits a decomposition
\begin{equation} \label{motive}
R(A/S) = \bigoplus_{i = 0}^{2g}R^i(A/S),
\end{equation}
with $[n]^*$ acting on $R^i(A/S)$ as multiplication by $n^i$. This gives the Beauville decomposition on the Chow ring ${\rm CH}(A)$: for $0 \leq i \leq g + d$, we have
\begin{equation}\label{beauville}
{\rm CH}^i(A) = \bigoplus_{j =  \max\{i - g, 2i - 2g\}}^{\min\{2i, i + d\}} {\rm CH}^{i}_{(j)}(A),
\end{equation}
where ${\rm CH}^i_{(j)}(A) := \big\{\alpha \in {\rm CH}^i(A) : [n]^*(\alpha) = n^{2i - j}\alpha, \textrm{ for all } n \in \mathbb{Z}\big\}$. Note that ${\rm CH}^i_{(j)}(A) = {\rm CH}^i\big(R^{2i - j}(A/S)\big)$.

Let $\mathcal{P}$ be a Poincaré line bundle on $A \times_S A$ rigidified along the zero sections, and denote its first Chern class by $\ell$. The Fourier transform $\mathcal{F} \colon {\rm CH}(A) \to {\rm CH}(A)$ is given by $\alpha \mapsto {\rm pr}_{2, *}\big({\rm pr}_1^*(\alpha) \cdot \exp(\ell) \big)$, where ${\rm pr}_1, {\rm pr}_2 \colon A \times_S A \to A$ are the two projections. There are isomorphisms
\begin{equation*}
\mathcal{F} \colon R^i(A/S) \xrightarrow{\sim} R^{2g - i}(A/S), \quad \mathcal{F} \colon {\rm CH}^i_{(j)}(A) \xrightarrow{\sim} {\rm CH}^{g - i + j}_{(j)}(A).
\end{equation*}

Here we find it more convenient to adopt the following notation:
\begin{equation}
{\rm CH}_{(i, j)}(A) := {\rm CH}^{\frac{i + j}{2}}_{(j)}(A),
\end{equation}
for all $i$ and $j$ such that $i + j$ is even and that the right-hand side makes sense. With this new notation, the first index $i$ reflects the weight in the motivic decomposition, and the second index $j$ reflects the level in the Beauville decomposition. In other words, we have
\begin{equation*}
{\rm CH}_{(i, j)}(A) = {\rm CH}\big(R^i(A/S)\big) \cap {\rm CH}_{(j)}(A).
\end{equation*}
Expression for the intersection (`$\cdot$') and Pontryagin (`$*$') products are even simpler: for $\alpha \in {\rm CH}_{(i, j)}(A)$ and $\beta \in {\rm CH}_{(k, l)}(A)$, we have $\alpha \cdot \beta \in {\rm CH}_{(i + k, j + l)}(A)$ and $\alpha * \beta \in {\rm CH}_{(i + k - g, j + l)}(A)$. Further, there are isomorphisms $\mathcal{F} \colon {\rm CH}_{(i, j)}(A) \xrightarrow{\sim} {\rm CH}_{(2g - i, j)}(A)$.

\subsection{\it The Dutch house} \label{dutchhouse} --- We introduce a useful picture to represent ${\rm CH}(A)$, called the {\it Dutch house}. It has the advantage of combining the Beauville decomposition, the motivic decomposition and the Fourier transform into one picture. It also enables us to make clear statements without complicated indices.

\begin{figure}
\centering
\includegraphics[height=.6\textheight]{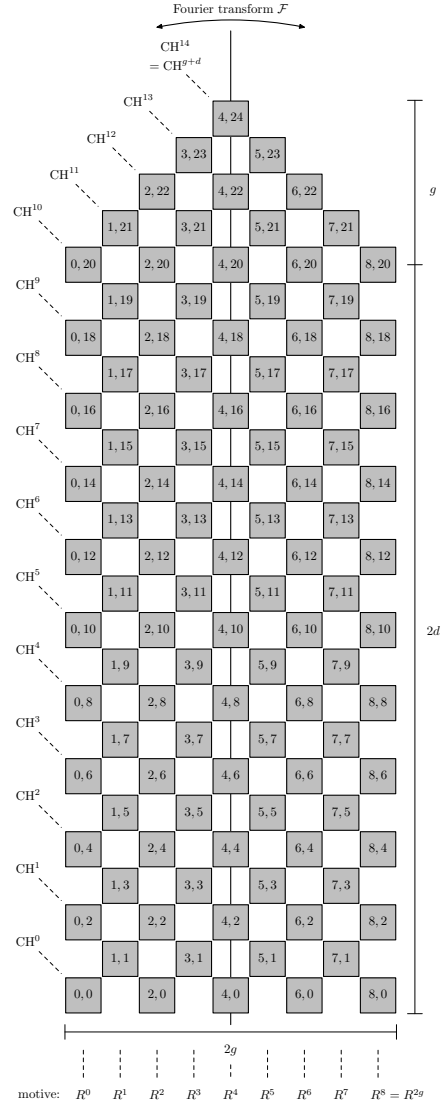}
\caption{The outside of the Dutch house ($g = 4$, $d = \dim\mathcal{M}_{4, 1} = 10$).}\label{outside}
\end{figure}

Figure~\ref{outside} is drawn based on the universal Jacobian over $\mathcal{M}_{4, 1}$. The abelian scheme $A$ has relative dimension $g = 4$, while the base $S$ is of dimension $d = \dim\mathcal{M}_{4, 1} = 10$. The $(i, j)$-th block represents the component ${\rm CH}_{(i, j)}(A)$ in ${\rm CH}(A)$. Therefore, the columns read the motivic decomposition, and the rows read levels in the Beauville decomposition. As a result of our indexing, components with the same codimension lie on a dashed line from top left to bottom right. Finally $\mathcal{F}$ acts as the symmetry with respect to the middle vertical line.

It is not difficult to verify that the house shape results from the precise index range of (\ref{beauville}). The width of the house depends on the relative dimension $g$, while the 
height (without roof) depends on the dimension of the base $d$. In particular when $S = {\rm Spec}(k)$, {\it i.e.}~when $d = 0$, the house degenerates to a pyramid ({\it cf.}~\cite{Moo09}, Figure~1). Note that {\it a priori\/} there are also components ${\rm CH}_{(i, j)}(A)$ with negative $j$. These components are expected to vanish according to the Beauville conjecture ({\it cf.}~\cite{Bea83}, Section 5). We have not drawn this negative part since the classes we shall study are all in ${\rm CH}_{(i, j)}(A)$ with $j \geq 0$, {\it i.e.}~inside the house.

\subsection{} For all $n \in \mathbb{Z}$, we have $[n]^* \pi^* = \pi^*$ and $\pi_* [n]_* = \pi_*$. This implies that the image of $\pi^* \colon {\rm CH}(S) \to {\rm CH}(A)$ is contained in $\oplus_{i = 0}^{d}{\rm CH}_{(0, 2i)}(A) = {\rm CH}\big(R^0(A/S)\big)$, and only elements of $\oplus_{i = 0}^{d}{\rm CH}_{(2g, 2i)}(A) = {\rm CH}\big(R^{2g}(A/S)\big)$ can have non-zero image under $\pi_* \colon {\rm CH}(A) \to {\rm CH}(S)$.

In fact we have even stronger results: let $\sigma_0 \colon S \to A$ be the zero section of the abelian scheme. Then there are isomorphisms of relative Chow motives that form a commutative diagram ({\it cf.}~\cite{DM91}, Example~1.4).
\begin{equation*}
\begin{tikzcd}[column sep=tiny]
\phantom{[2g]} R^0(A/S) \arrow[leftarrow]{rr}{\mathcal{F}}[swap]{\sim} & & R^{2g}(A/S)[2g] \arrow{dl}{\sigma_0^*}[swap]{\sim} \\
& R(S/S) \arrow{ul}{\pi^*}[swap]{\sim}
\end{tikzcd}
\end{equation*}
At the level of Chow groups, we obtain a commutative diagram of isomorphisms of $\mathbb{Q}$-algebras.
\begin{equation} \label{isom}
\begin{tikzcd}[column sep=tiny]
\big(\bigoplus_{i = 0}^{d}{\rm CH}_{(0, 2i)}(A), \cdot\big) \arrow[leftarrow]{rr}{\mathcal{F}}[swap]{\sim} & & \big(\bigoplus_{i = 0}^{d}{\rm CH}_{(2g, 2i)}(A), *\big) \arrow{dl}{\pi_*}[swap]{\sim} \\
& \big({\rm CH}(S), \cdot\big) \arrow{ul}{\pi^*}[swap]{\sim}
\end{tikzcd}
\end{equation}
The gradings are preserved as $\pi^* \colon {\rm CH}^i(S) \xrightarrow{\sim} {\rm CH}_{(0, 2i)}(A)$ and $\pi_* \colon {\rm CH}_{(2g, 2i)}(A) \xrightarrow{\sim} {\rm CH}^i(S)$.

In particular, the Chow ring ${\rm CH}(S)$ may be regarded as a $\mathbb{Q}$-subalgebra of $\big({\rm CH}(A), \cdot\big)$ via $\pi^*$, or as a $\mathbb{Q}$-subalgebra of $\big({\rm CH}(A), *\big)$ via $\pi_*$. In terms of the Dutch house, we may identify ${\rm CH}(S)$ with the $0$-th column or with the $2g$-th column of the house.

\subsection{\it The \texorpdfstring{$\mathfrak{sl}_2$}{sl2}-action} --- We refer to Künnemann's work \cite{Kün93} for the motivic Lefschetz decomposition of a polarized abelian scheme. The Lefschetz and Lambda operators described in {\it loc.\,cit.}~generate an $\mathfrak{sl}_2$-action on the relative Chow motive and on the Chow ring of the abelian scheme. 

From now on we shall consider the following setting. Take $S$ a smooth connected variety of dimension $d$ over $k$. Let $p \colon C \to S$ be a smooth relative curve of genus $g \geq 2$ over $S$, together with a section $x_0 \colon S \to C$. Denote by $\pi \colon J \to S$ the associated  relative Jacobian, which is an abelian scheme of relative dimension $g$. Let $\iota \colon C \hookrightarrow J$ be the closed embedding given on points by $x \mapsto \mathcal{O}_C(x - x_0)$. Then if we write $\sigma_0 \colon S \to J$ for the zero section of $J$, we have $\iota \circ x_0 = \sigma_0$. To summarize, we shall work with the commutative diagram below.
\begin{equation} \label{setting}
\begin{tikzcd}[row sep=large, column sep=scriptsize]
C \arrow[hookrightarrow]{rr}{\iota} \arrow{dr}{p} & & J \arrow{dl}[swap]{\pi} \\
& S \arrow[bend left=50]{ul}{x_0} \arrow[bend right=50]{ur}[swap]{\sigma_0}
\end{tikzcd}
\end{equation}
We write $[x_0] := \big[x_0(S)\big] \in {\rm CH}^1(C)$, $[\sigma_0] := \big[\sigma_0(S)\big] \in {\rm CH}^g(J)$ and $[C] := \big[\iota(C)\big] \in {\rm CH}^{g - 1}(J)$ for short.

In \cite{Pol07b}, Polishchuk showed that the $\mathfrak{sl}_2$-action can be constructed geometrically using the curve class $[C] \in {\rm CH}^{g - 1}(J)$. By translating his notation into ours, we have the following: consider the decomposition
\begin{equation*}
[C] = \sum_{j = 0}^{2g - 2} [C]_{(j)},
\end{equation*}
with $[C]_{(j)} \in {\rm CH}^{g - 1}_{(j)}(J) = {\rm CH}_{(2g - 2 - j, j)}(J)$. We take the component $[C]_{(0)}$ and write 
\begin{equation}
\theta := - \mathcal{F}\big([C]_{(0)}\big) \in {\rm CH}^1_{(0)}(J) = {\rm CH}_{(2, 0)}(J).
\end{equation}

\subsection{\it Proposition {\rm (Polishchuk)}} --- {\it The divisor class $\theta$ is relatively ample and induces the natural principal polarization of $J$. Define operators $e, f$ and $h$ on ${\rm CH}(J)$ by}
\begin{align*}
e(\alpha) & := -\theta \cdot \alpha, \\
f(\alpha) & := -[C]_{(0)} * \alpha, \\
h(\alpha) & := (i - g) \cdot \alpha, \textrm{ for } \alpha \in {\rm CH}_{(i, j)}(J).
\end{align*}
{\it Then the operators above generate a $\mathbb{Q}$-linear representation $\mathfrak{sl}_2 \to {\rm End}_{\mathbb{Q}}\big({\rm CH}(J)\big)$.}

\medskip
In terms of the Dutch house, the operator $e$ shifts classes to the right by $2$ blocks, while $f$ shifts classes to the left by $2$ blocks. The middle column of the house has weight $0$ with respect to $h$.

We refer to \cite{Pol07b}, Theorem~2.6 for the proof. Note that $\theta$ is equal to $-\tau_2(C)/2$ in Polishchuk's notation. Also we followed Polishchuk's convention by defining $e = -\theta$, which differs by a sign from Künnemann's version ({\it cf.}~Remark after Theorem~2.6 in {\it loc.\,cit.}).

Moreover, the Fourier transform $\mathcal{F}$ can be written as
\begin{equation} \label{fourier}
\mathcal{F} = \exp(e) \exp(-f) \exp(e).
\end{equation}
This is proven via the same argument as in \cite{Bea04},~2.7, using the identity
\begin{equation*}
\ell = {\rm pr}_1^*(\theta) + {\rm pr}_2^*(\theta) - m^*(\theta).
\end{equation*}
See also \cite{Pol07}, Section~1, where Polishchuk explained that $\mathcal{F}$ corresponds to the action of the matrix $\big(\begin{smallmatrix} 0 & 1 \\ -1 & 0 \end{smallmatrix}\big)$ in the algebraic group ${\rm SL}_2$.

The following lemma shows that $\theta$ is the class of a relative theta divisor in the classical sense: it comes from a family of theta characteristics (with $\mathbb{Q}$-coefficients).

\subsection{\it Lemma} \label{pullback} --- {\it We have $\iota^*(\theta) = K/2 + [x_0] + \psi/2$ in ${\rm CH}^1(C)$.}

\begin{proof}
The goal is to calculate $\iota^*(\theta) = -\iota^*\Big(\mathcal{F}\big([C]_{(0)}\big)\Big)$ and we start from $i^*\Big(\mathcal{F}\big([C]\big)\Big)$. Consider the following three Cartesian squares
\begin{equation*}
\begin{tikzcd}[row sep=large, column sep=large]
C \times_S C \arrow{r}{(\iota, {\rm id}_C)} \arrow{d}[swap]{({\rm id}_C, \iota)} & J \times_S C \arrow{r}{{\rm pr}_2} \arrow{d}{({\rm id}_J, \iota)} & C \arrow{d}{\iota} \\
C \times_S J \arrow{r}[swap]{(\iota, {\rm id}_J)} \arrow{d}[swap]{{\rm pr}_1} & J \times_S J \arrow{r}[swap]{{\rm pr}_2} \arrow{d}{{\rm pr}_1} & J \\
C \arrow{r}[swap]{\iota} & J  &  \phantom{C \times_S C}
\end{tikzcd}
\end{equation*}
where ${\rm pr}_1$ and ${\rm pr}_2$ stand for the two projections in all cases. Then we have
\begin{align*}
\iota^*\Big(\mathcal{F}\big([C]\big)\Big) & = \iota^*{\rm pr}_{2, *}\big({\rm pr}_1^*\,\iota_*(\mathbbm{1}_C) \cdot \exp(\ell)\big) \\
& = {\rm pr}_{2, *}({\rm id}_J, \iota)^*\big({\rm pr}_1^*\,\iota_*(\mathbbm{1}_C) \cdot \exp(\ell)\big) \tag{Cartesian square on the right} \\
& = {\rm pr}_{2, *}({\rm id}_J, \iota)^*\big((\iota, {\rm id}_J)_*{\rm pr}_1^*(\mathbbm{1}_C) \cdot \exp(\ell)\big) \tag{Cartesian square at the bottom} \\
& = {\rm pr}_{2, *}({\rm id}_J, \iota)^*\big((\iota, {\rm id}_J)_*(\mathbbm{1}_{C \times_S J}) \cdot \exp(\ell)\big) \nonumber \\
& = {\rm pr}_{2, *}({\rm id}_J, \iota)^*(\iota, {\rm id}_J)_*(\iota, {\rm id}_J)^*\big(\exp(\ell)\big) \nonumber \\
& = {\rm pr}_{2, *}(\iota, {\rm id}_C)_*({\rm id}_C, \iota)^*(\iota, {\rm id}_J)^*\big(\exp(\ell)\big) \tag{Cartesian square on the top-left} \\
& = {\rm pr}_{2, *}(\iota, \iota)^*\big(\exp(\ell)\big) \nonumber \\
& = {\rm pr}_{2, *}\Big(\exp\big((\iota, \iota)^*(\ell)\big)\Big). \nonumber
\end{align*}
The class $\ell$ has the property that
\begin{equation*}
(\iota, \iota)^*(\ell) = \Delta - {\rm pr}_1^*\big([x_0]\big) - {\rm pr}_2^*\big([x_0]\big),
\end{equation*}
with $\Delta$ the class of the diagonal in $C \times_S C$. It follows that
\begin{align} \label{expand}
\iota^*\Big(\mathcal{F}\big([C]\big)\Big) & = {\rm pr}_{2, *}\bigg(\exp\Big(\Delta - {\rm pr}_1^*\big([x_0]\big) - {\rm pr}_2^*\big([x_0]\big)\Big)\bigg) \\
& = {\rm pr}_{2, *}\bigg(\exp\Big(\Delta - {\rm pr}_1^*\big([x_0]\big)\Big)\bigg) \cdot \exp\big(\!-\![x_0]\big). \nonumber
\end{align}
Then we observe that on the left-hand side of (\ref{expand}), we have
\begin{equation*}
\iota^*\Big(\mathcal{F}\big([C]\big)\Big) = \sum_{j = 0}^{2g - 2} \iota^*\Big(\mathcal{F}\big([C]_{(j)}\big)\Big),
\end{equation*}
with $\iota^*\Big(\mathcal{F}\big([C]_{(j)}\big)\Big) \in {\rm CH}^{j + 1}(C)$. Hence $\iota^*\Big(\mathcal{F}\big([C]_{(0)}\big)\Big)$ is the component of $\iota^*\Big(\mathcal{F}\big([C]\big)\Big)$ in codimension $1$. Expanding the right-hand side of (\ref{expand}), we find
\begin{align*}
\iota^*\Big(\mathcal{F}\big([C]_{(0)}\big)\Big) & = {\rm pr}_{2, *}\bigg(\frac{1}{2}\Big(\Delta - {\rm pr}_1^*\big([x_0]\big)\Big)^2\bigg) - {\rm pr}_{2, *}\Big(\Delta - {\rm pr}_1^*\big([x_0]\big)\Big) \cdot [x_0] \\
& = \frac{1}{2}{\rm pr}_{2, *}(\Delta \cdot \Delta) - {\rm pr}_{2, *}\Big(\Delta \cdot {\rm pr}_1^*\big([x_0]\big)\Big) + \frac{1}{2}{\rm pr}_{2, *}{\rm pr}_1^*\big([x_0] \cdot [x_0]\big) \nonumber \\
& \qquad - {\rm pr}_{2, *}(\Delta) \cdot [x_0] + {\rm pr}_{2, *}{\rm pr}_1^*\big([x_0]\big) \cdot [x_0]. \nonumber
\end{align*}
Four of the five terms in the previous expression are easily calculated:
\begin{gather*}
{\rm pr}_{2, *}(\Delta \cdot \Delta) = -K, \\
{\rm pr}_{2, *}\Big(\Delta \cdot {\rm pr}_1^*\big([x_0]\big)\Big) = {\rm pr}_{2, *}(\Delta) \cdot [x_0] = {\rm pr}_{2, *}{\rm pr}_1^*\big([x_0]\big) \cdot [x_0] = [x_0].
\end{gather*}
For the term ${\rm pr}_{2, *}{\rm pr}_1^*\big([x_0] \cdot [x_0]\big)$, consider the following Cartesian square.
\begin{equation*}
\begin{tikzcd}[row sep=large, column sep=large]
C \times_S C \arrow{r}{{\rm pr}_2} \arrow{d}[swap]{{\rm pr}_1} & C \arrow{d}{p} \\
C \arrow{r}[swap]{p} & S
\end{tikzcd}
\end{equation*}
Then we have
\begin{equation*}
{\rm pr}_{2, *}{\rm pr}_1^*\big([x_0] \cdot [x_0]\big) = p^*p_*\big([x_0] \cdot [x_0]\big) = p^*p_*x_{0, *}x_0^*\big([x_0]\big) = p^*x_0^*\big([x_0]\big) = \psi.
\end{equation*}
In total $\iota^*\Big(\mathcal{F}\big([C]_{(0)}\big)\Big) = -K/2 - [x_0] - \psi/2 - [x_0] + [x_0] = -K/2 - [x_0] - \psi/2$.
\end{proof}

\medskip
\section{\bf The tautological ring of a relative Jacobian}
\medskip

\noindent Throughout this section, we work in the setting of diagram~(\ref{setting}). We shall define the tautological ring $\mathcal{T}(J)$ of the relative Jacobian $J$ and describe its structure. The central result is that via the pull-back map $\pi^*$, the tautological ring $\mathcal{R}(S)$ can be identified with the $\mathbb{Q}$-subalgebra of $\mathcal{T}(J)$ located on the $0$-th column of the Dutch house. At the end of this section, we give an algebraic proof of an identity obtained by Morita.

\subsection{\it Definition} \label{defTJ} --- We define the {\it tautological ring\/} of $J$ as the smallest $\mathbb{Q}$-subalgebra of $\big({\rm CH}(J), \cdot\big)$ that contains the class $[C] \in {\rm CH}^{g - 1}(J)$, and that is stable under the Fourier transform $\mathcal{F}$ and the action of $[n]^*$. We denote it by $\mathcal{T}(J)$, and we write $\mathcal{T}_{(i, j)}(J) := \mathcal{T}(J) \cap {\rm CH}_{(i, j)}(J)$.

\medskip
Note that the notion of `tautological ring' was first invented by Beauville \cite{Bea04}, working on a Jacobian variety and modulo algebraic equivalence. Since then there have been various versions of the tautological ring ({\it cf.}~\cite{Pol07} and \cite{Moo09}). In \cite{Pol07b}, Polishchuk considered in the relative setting a much bigger tautological ring, including all classes of $\pi^*\big({\rm CH}(S)\big)$. Here our definition works best for the purpose of studying the tautological ring of $S$.

\subsection{} \label{pij} We make a few remarks on the definition above. First, the fact that $\mathcal{T}(J)$ is stable under `$\cdot$' and $\mathcal{F}$ implies that it is also stable under the Pontryagin product `$*$'. Second, saying that $\mathcal{T}(J)$ is stable under the action of $[n]^*$ is equivalent to saying that it is stable under decomposition~(\ref{beauville}).

\begin{figure}
\centering
\includegraphics[height=.43\textheight]{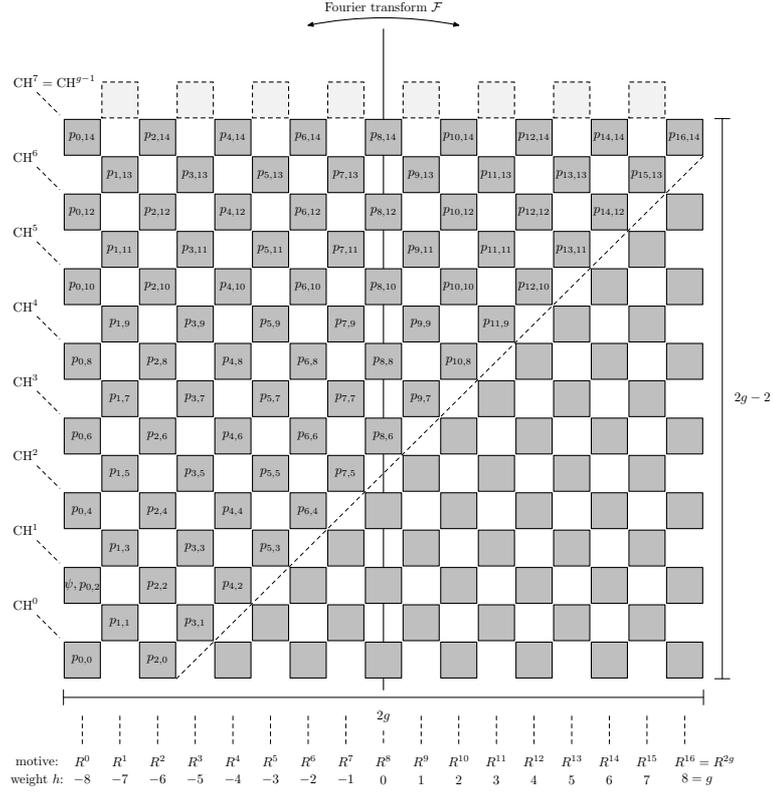}
\caption{The inside of the Dutch house ($g = 8$).}
\label{inside}
\end{figure}

In particular, we have $[C]_{(j)} \in \mathcal{T}_{(2g - 2 - j, j)}(J)$, which then implies that $\theta = -\mathcal{F}\big([C]_{(0)}\big) \in \mathcal{T}_{(2, 0)}(J)$. For $0 \leq j \leq 2g - 2$ and $0 \leq i \leq j + 2$ such that $i + j$ is even, consider the class
\begin{equation*}
\theta^{(j - i + 2)/2} \cdot [C]_{(j)} \in \mathcal{T}_{(2g - i, j)}(J).
\end{equation*}
Denote its Fourier dual by
\begin{equation}
p_{i, j} := \mathcal{F}\big(\theta^{(j - i + 2)/2} \cdot [C]_{(j)}\big) \in \mathcal{T}_{(i, j)}(J).
\end{equation}
As examples we have $p_{2, 0} = \mathcal{F}\big([C]_{(0)}\big) = -\theta$ and $p_{0, 0} = \mathcal{F}\big(\theta \cdot [C]_{(0)}\big) = g \cdot \mathbbm{1}_J$. In particular, the action of $e \in \mathfrak{sl}_2$ on ${\rm CH}(J)$ is just the intersection with $p_{2, 0}$. We extend the definition naturally by setting $p_{i, j} := 0$ for $i < 0$ or $j < 0$ or $j > 2g - 2$ (again $i \leq j + 2$ and $i + j$ even). It is not difficult to see that\begin{equation*}
f(p_{i, j}) = p_{i - 2, j}.
\end{equation*}
Figure~\ref{inside} depicts the classes $\{p_{i, j}\}$ inside the Dutch house for $g = 8$, together with the class $\psi \in {\rm CH}_{(0, 2)}(J)$. Note that when $d = \dim(S)$ is small, those classes that are `above the roof' vanish for obvious reasons.

\subsection{\it Theorem} \label{structure} --- (i) {\it The tautological ring $\mathcal{T}(J)$ coincides with the $\mathbb{Q}$-subalgebra of $\big({\rm CH}(J), \cdot\big)$ generated by the classes $\{p_{i, j}\}$ and $\psi$.}

\smallskip
\noindent\makebox[\leftmargin][r]{(ii) }{\it The operator $f \in \mathfrak{sl}_2$ acts on polynomials in $\{p_{i, j}\}$ and $\psi$ via the following differential operator of degree $2$:}
\begin{equation}
\begin{split}
\mathcal{D}  :=  \frac{1}{2}\sum_{\substack{i, j \\ k, l}}\Bigg(\psi p_{i - 1, j - 1} p_{k - 1, l - 1} & - \binom{i + k - 2}{i - 1} p_{i + k - 2, j + l}\Bigg) \partial p_{i, j} \partial p_{k, l} \\
& + \sum_{i, j}p_{i - 2, j}\partial p_{i, j}. 
\end{split}
\end{equation}

\medskip
We may track the differential operator $\mathcal{D}$ in the Dutch house. For the degree $2$ part of $\mathcal{D}$, whenever there is a product of two classes $p_{i, j}p_{k, l}$, first find the generators to the lower-left of $p_{i, j}$ and $p_{k, l}$ by $1$ block, and multiply by $\psi$, which yields $\psi p_{i - 1, j - 1}p_{k - 1, l - 1}$. Then look for the generator to the left of $p_{i, j}p_{k, l}$ by $2$ blocks, which is $p_{i + k - 2, j + l}$. For the linear part of $\mathcal{D}$, the generator $p_{i, j}$ is simply replaced by the one to the left of it by $2$ blocks, {\it i.e.}~$p_{i - 2, j}$. All these operations shift classes to the left by $2$ blocks.

\begin{proof}
Suppose that we have proven (ii) and that $\psi \in \mathcal{T}(J)$. Consider the $\mathbb{Q}$-subalgebra of $\big({\rm CH}(J), \cdot\big)$ generated by the classes $\{p_{i, j}\}$ and $\psi$. We denote it by $\mathcal{T}'(J)$ and we have $\mathcal{T}'(J) \subset \mathcal{T}(J)$. By definition $\mathcal{T}'(J)$ is stable under decomposition~(\ref{beauville}). It is stable under the action of $e \in \mathfrak{sl}_2$, which is the intersection with $p_{2, 0}$. Then by (ii), we know that $\mathcal{T}'(J)$ is also stable under the action of $f \in \mathfrak{sl}_2$. It follows that $\mathcal{T}'(J)$ is stable under the Fourier transform $\mathcal{F} = \exp(e)\exp(-f)\exp(e)$. In particular, the classes $\big\{[C]_{(j)}\big\}$ are contained in $\mathcal{T}'(J)$. Since $\mathcal{T}(J)$ is defined as the smallest $\mathbb{Q}$-algebra that satisfies these properties, there is necessarily an equality $\mathcal{T}'(J) = \mathcal{T}(J)$, which proves (i).

Statement (ii) follows essentially from \cite{Pol07b}, Formula~2.9. We just need to translate the notation carefully. Following Polishchuk, we write $\eta := K/2 + [x_0] + \psi/2$, which is equal to $\iota^*(\theta)$ by Lemma~\ref{pullback}. First of all, we have $f = \displaystyle-\frac{1}{2}\tilde{X}_{2, 0}(C)$ in his notation. We define operators $\tilde{p}_{i, j}$ on ${\rm CH}(J)$ by $\tilde{p}_{i, j}(\alpha) = p_{i, j} \cdot \alpha$. Then the fact that 
\begin{equation*}
p_{i, j} = \mathcal{F}\big(\theta^{(j - i + 2)/2} \cdot [C]_{(j)}\big) = \mathcal{F}\big(\iota_*(\eta^{(j - i + 2)/2})_{(j)}\big)
\end{equation*}
is translated into
\begin{equation} \label{translate}
\tilde{p}_{i, j} = \frac{1}{i!}\tilde{X}_{0, i}(\eta^{(j - i + 2)/2}).
\end{equation}
We apply Formula~2.9 in {\it loc.\,cit.}~and find
\begin{align} \label{step1}
[f, \tilde{p}_{i, j}] & = -\frac{1}{2 \cdot i!}\big[\tilde{X}_{2, 0}(C), \tilde{X}_{0, i}(\eta^{(j - i + 2)/2})\big] \\
& = \frac{1}{(i - 1)!}\tilde{X}_{1, i - 1}(\eta^{(j - i + 2)/2}) - \frac{1}{(i - 2)!}\tilde{X}_{0, i - 2}(\eta^{(j - i + 4)/2}). \nonumber
\end{align}
Note that the second equality of (\ref{step1}) also involves the fact that $\tilde{X}_{i, 0}(C) = 0$ for $i \leq 1$ (Lemma~2.8 in {\it loc.\,cit.}), and that $x_0^*(\eta) = x_0^*\,\iota^*(\theta) = \sigma_0^*(\theta) = 0$ (since $\theta \in {\rm CH}_{(2, 0)}(J)$). We continue to calculate
\begin{align*}
\big[[f, \tilde{p}_{i, j}], \tilde{p}_{k, l}\big] & = \frac{1}{(i - 1)!\,k!}\big[\tilde{X}_{1, i - 1}(\eta^{(j - i + 2)/2}), \tilde{X}_{0, k}(\eta^{(l - k + 2)/2})\big] \\
& \qquad - \frac{1}{(i - 2)!\,k!}\big[\tilde{X}_{0, i - 2}(\eta^{(j - i + 4)/2}), \tilde{X}_{0, k}(\eta^{(l - k + 2)/2})\big]. \nonumber
\end{align*}
By applying the same formula, we have $\big[\tilde{X}_{0, i - 2}(\eta^{(j - i + 4)/2}), \tilde{X}_{0, k}(\eta^{(l - k + 2)/2})\big] = 0$, and
\begin{align*}
\big[\tilde{X}_{1, i - 1}(\eta^{(j - i + 2)/2}), \tilde{X}_{0, k}(\eta^{(l - k + 2)/2})\big] & = k \psi \tilde{X}_{0, k - 1}(\eta^{(l - k + 2)/2}) \tilde{X}_{0, i - 1}(\eta^{(j - i + 2)/2}) \\
& \qquad -k \tilde{X}_{0, i + k - 2}(\eta^{(j - i + l - k + 4)/2}). \nonumber
\end{align*}
In total, we obtain
\begin{align} \label{degree2}
\big[[f, \tilde{p}_{i, j}], \tilde{p}_{k, l}\big] & = \frac{1}{(i - 1)!\,(k - 1)!}\Big(\psi \tilde{X}_{0, k - 1}(\eta^{(l - k + 2)/2}) \tilde{X}_{0, i - 1}(\eta^{(j - i + 2)/2}) \\
& \qquad - \tilde{X}_{0, i + k - 2}(\eta^{(j - i + l - k + 4)/2})\Big) \nonumber \\
& = \psi \tilde{p}_{k - 1, l - 1} \tilde{p}_{i - 1, j - 1} - \binom{i + k - 2}{i - 1} \tilde{p}_{i + k - 2, j + l} \tag{by (\ref{translate})} \\
& = \psi \tilde{p}_{i - 1, j - 1} \tilde{p}_{k - 1, l - 1} - \binom{i + k - 2}{i - 1} \tilde{p}_{i + k - 2, j + l}. \tag{$\tilde{p}_{i - 1, j - 1}$ and $\tilde{p}_{k - 1, l - 1}$ commute}
\end{align}
On the other hand, since $f(\mathbbm{1}_J) = 0$, we have
\begin{equation} \label{degree1}
[f, \tilde{p}_{i, j}](\mathbbm{1}_J) = f(p_{i, j}) = p_{i - 2, j}.
\end{equation}
The relations (\ref{degree2}) and (\ref{degree1}) then imply that for any polynomial $P$ in $\{p_{i, j}\}$ and $\psi$, we have
\begin{equation*}
f\Big(P\big(\{p_{i, j}\}, \psi\big)\Big) = \mathcal{D}\Big(P\big(\{p_{i, j}\}, \psi\big)\Big),
\end{equation*}
where $\mathcal{D}$ is the differential operator stated in the theorem ({\it cf.}~\cite{Pol07}, Section~3).

It remains to prove that $\psi \in \mathcal{T}(J)$. To see this, we apply $\mathcal{D}$ to the class $p_{1, 1}^2 \in \mathcal{T}(J)$:
\begin{equation*}
\mathcal{D}(p_{1, 1}^2) = \psi p_{0, 0}^2 - \binom{0}{0} p_{0, 2} = g^2\psi - p_{0, 2} \in \mathcal{T}(J).
\end{equation*}
Hence $\psi = \big(\mathcal{D}(p_{1, 1}^2) + p_{0, 2}\big)/g^2 \in \mathcal{T}(J)$.
\end{proof}

\subsection{\it Corollary} \label{subring}--- {\it For $0 \leq i \leq g - 1$, there is an identity}
\begin{equation} \label{basechange}
p_{0, 2i} = \pi^*\Bigg(\frac{1}{2^{i + 1}}\sum_{0 \leq j \leq i} \binom{i + 1}{j + 1} \psi^{i - j} \kappa_j + \psi^i\Bigg).
\end{equation}
{\it Moreover, we have the following isomorphisms of $\mathbb{Q}$-algebras.}
\begin{equation} \label{tautisom}
\begin{tikzcd}[column sep=tiny]
\big(\bigoplus_{i = 0}^{d}\mathcal{T}_{(0, 2i)}(J), \cdot\big) \arrow[leftarrow]{rr}{\mathcal{F}}[swap]{\sim} & & \big(\bigoplus_{i = 0}^{d}\mathcal{T}_{(2g, 2i)}(J), *\big) \arrow{dl}{\pi_*}[swap]{\sim} \\
& \big(\mathcal{R}(S), \cdot\big) \arrow{ul}{\pi^*}[swap]{\sim}
\end{tikzcd}
\end{equation}
{\it In particular, the tautological ring $\mathcal{R}(S)$ may be regarded as a $\mathbb{Q}$-subalgebra of $\big(\mathcal{T}(J), \cdot\big)$ via $\pi^*$.}

\begin{proof}
By (\ref{isom}) we have $p_{0, 2i} = \mathcal{F}\big(\theta^{i + 1} \cdot [C]_{(2i)}\big) = \pi^*\pi_*\big(\theta^{i + 1} \cdot [C]_{(2i)}\big) = \pi^*\pi_*\big(\theta^{i + 1} \cdot [C]\big)$, hence it suffices to calculate $\pi_*\big(\theta^{i + 1} \cdot [C]\big)$. Then Lemma~\ref{pullback} and the projection formula imply that
\begin{align*}
\pi_*\big(\theta^{i + 1} \cdot [C]\big) & = p_*\Bigg(\bigg(\frac{1}{2}K + [x_0] + \frac{1}{2}\psi\bigg)^{i + 1}\Bigg) \\
& = \sum_{\substack{j + k + l = i + 1 \\ j, k, l \geq 0}} \frac{(i + 1)!}{j!\,k!\,l!}\frac{1}{2^{j + l}} \, p_*\big(K^j \cdot [x_0]^k \cdot \psi^l\big). \nonumber
\end{align*}
Again by applying the projection formula to $p \colon C \to S$ and $x_0 \colon S \to C$, we find
\begin{equation}
p_*\big(K^j \cdot [x_0]^k \cdot \psi^l\big) = \psi^l \cdot p_*\big(K^j \cdot [x_0]^k\big) = 
\begin{cases}
\psi^l \cdot \kappa_{j - 1} & \textrm{if } k = 0, \\
\psi^l \cdot x_0^*\big(K^j \cdot [x_0]^{k - 1}\big) = (-1)^{k - 1} \psi^i & \textrm{if } k \geq 1,
\end{cases}
\end{equation}
with the convention $\kappa_{-1} = 0$. It follows that
\begin{align*}
\pi_*\big(\theta^{i + 1} \cdot [C]\big) & = \sum_{\substack{j + l = i + 1 \\ j, l \geq 0}} \frac{(i + 1)!}{j!\,l!}\frac{1}{2^{i + 1}} \, \psi^l \cdot \kappa_{j - 1} + \sum_{\substack{j + k + l = i + 1 \\ j, k, l \geq 0}} \frac{(i + 1)!}{j!\,k!\,l!}\frac{1}{2^{j + l}} (-1)^{k - 1} \psi^i \\
& \qquad - \sum_{\substack{j + l = i + 1 \\ j, l \geq 0}} \frac{(i + 1)!}{j!\,l!}\frac{1}{2^{i + 1}} (-1) \psi^i \nonumber \\
& = \frac{1}{2^{i + 1}}\sum_{0 \leq j \leq i} \binom{i + 1}{j + 1} \psi^{i - j} \kappa_j + \bigg(\frac{1}{2} - 1 + \frac{1}{2}\bigg)^{i + 1} \psi^i + \bigg(\frac{1}{2} + \frac{1}{2}\bigg)^{i + 1} \psi^i \nonumber \\
& = \frac{1}{2^{i + 1}}\sum_{0 \leq j \leq i} \binom{i + 1}{j + 1} \psi^{i - j} \kappa_j + \psi^i, \nonumber
\end{align*}
which proves the identity (\ref{basechange}).

In particular, since $\oplus_{i = 0}^{d}\mathcal{T}_{(0, 2i)}(J)$ is generated by the classes $\{p_{0, 2i}\}$ and $\psi$, we have an inclusion $\oplus_{i = 0}^{d}\mathcal{T}_{(0, 2i)}(J) \subset \pi^*\big(\mathcal{R}(S)\big)$. To see the other inclusion, we know from (\ref{basechange}) that for $0 \leq j \leq g - 1$,  the class $\pi^*(\kappa_j)$ can also be expressed as linear combinations of $\{p_{0, 2i}\}$ and $\psi$. For $j > g - 1$, there is a classical result of Mumford \cite{Mum83} saying that $\kappa_j = 0$ for $j \geq g - 1$. Therefore we have $\oplus_{i = 0}^{d}\mathcal{T}_{(0, 2i)}(J) = \pi^*\big(\mathcal{R}(S)\big)$, and the rest follows from (\ref{isom}).
\end{proof}

\subsection{\it Example: a theorem of Morita} --- We conclude this section by proving an identity of Morita. The proof has the advantage of being purely algebraic, which holds over fields of arbitrary characteristic. But as everything else in this note, it only works with $\mathbb{Q}$-coefficients.

The identity reveals some connection between $\mathcal{T}^{g + 1}(J)$ and $\mathcal{R}^1(S)$. Consider the class $[C]_{(1)} \in \mathcal{T}_{(2g - 3, 1)}(J)$ and its Fourier dual $\mathcal{F}\big([C]_{(1)}\big) \in \mathcal{T}_{(3, 1)}(J)$. Then we have $[C]_{(1)} \cdot \mathcal{F}\big([C]_{(1)}\big) \in \mathcal{T}_{(2g, 2)}(J)$, and Morita answered what the image under $\pi_*$ of this class is. We refer to the paper of Hain-Reed \cite{HR01} for the original statement, which (with $\mathbb{Q}$-coefficient) is equivalent to the following.

\subsection{\it Theorem {\rm (Morita)}} \label{morita} --- {\it We have $\pi_*\Big([C]_{(1)} \cdot \mathcal{F}\big([C]_{(1)}\big)\Big) = \kappa_1/6 + g\psi$ in $\mathcal{R}^1(S)$.}

\medskip
Note that $\kappa_1/12$ is equal to $\lambda_1$, which stands for the first Chern class of the Hodge bundle $p_*(\Omega_{C/S})$. Hence the right-hand side is also equal to $2\lambda_1 + g\psi$.

\begin{proof}
Recall that $\mathcal{F}\big([C]_{(1)}\big) = p_{3, 1}$ and that by (\ref{isom}), we have
\begin{equation} \label{goal}
\pi^*\pi_*\big([C]_{(1)} \cdot p_{3, 1}\big) = \mathcal{F}\big([C]_{(1)} \cdot p_{3, 1}\big) = -p_{3, 1} * [C]_{(1)}.
\end{equation}
Therefore it suffices to express $-p_{3, 1} * [C]_{(1)}$ in terms of $p_{0, 2}$ and $\psi$, and then apply Corollary~\ref{subring}. 

The first step is to express $[C]_{(1)}$ in terms of the classes $\{p_{i, j}\}$. By definition we have
\begin{equation*}
f(p_{3, 1}) = p_{1, 1}, \quad fe(p_{1, 1}) = ef(p_{1, 1}) - h(p_{1, 1}) = (g - 1)p_{1, 1}.
\end{equation*}
Therefore $f\Big(p_{3, 1} - \frac{1}{g - 1} e(p_{1, 1})\Big) = 0$, which implies 
\begin{equation*}
e^{g - 2}\bigg(p_{3, 1} - \frac{1}{g - 1} e(p_{1, 1})\bigg) = e^{g - 2}(p_{3, 1}) - \frac{1}{g - 1}e^{g - 1}(p_{1, 1}) = 0.
\end{equation*}
Apply $\mathcal{F}^{-1}$ to the previous equation, and we find
\begin{equation*}
(-1)^{g - 2}f^{g - 2}\big([C]_{(1)}\big) - \frac{(-1)^{g - 1}}{g - 1}f^{g - 1}(-e)\big([C]_{(1)}\big) = 0,
\end{equation*}
so that
\begin{equation} \label{irred}
f^{g - 2}\big([C]_{(1)}\big) - \frac{1}{g - 1}f^{g - 1}e\big([C]_{(1)}\big) = 0.
\end{equation}
On the other hand, by (\ref{fourier}) we have
\begin{align*}
p_{3, 1}  =  \mathcal{F}\big([C]_{(1)}\big) & =  \exp(e)\exp(-f)\exp(e)\big([C]_{(1)}\big) \\
& = \exp(e)\exp(-f)\Big([C]_{(1)} + e\big([C]_{(1)}\big)\Big) \nonumber \\
& = \exp(e)\Bigg(\frac{(-1)^{g - 3}}{(g - 3)!} f^{g - 3}\big([C]_{(1)}\big) + \frac{(-1)^{g - 2}}{(g - 2)!} f^{g - 2}\big([C]_{(1)}\big) \nonumber \\
& \qquad + \frac{(-1)^{g - 2}}{(g - 2)!} f^{g - 2}e\big([C]_{(1)}\big) + \frac{(-1)^{g - 1}}{(g - 1)!} f^{g - 1}e\big([C]_{(1)}\big)\Bigg) \nonumber \\
& = \exp(e)\bigg(\frac{(-1)^{g - 3}}{(g - 3)!} f^{g - 3}\big([C]_{(1)}\big) + \frac{(-1)^{g - 2}}{(g - 2)!} f^{g - 2}e\big([C]_{(1)}\big)\bigg) \tag{by (\ref{irred})} \\
& = \frac{(-1)^{g - 3}}{(g - 3)!} f^{g - 3}\big([C]_{(1)}\big) + \frac{(-1)^{g - 2}}{(g - 2)!} f^{g - 2}e\big([C]_{(1)}\big).\nonumber
\end{align*}
Note that for $g = 2$, we ignore the $f^{g - 3}$ term and the rest of the argument still works. Then we apply $\mathcal{F}$ to both sides, which gives
\begin{align*}
(-1)^{g + 1}[C]_{(1)} & = \frac{1}{(g - 3)!}e^{g - 3}(p_{3, 1}) - \frac{1}{(g - 2)!}e^{g - 2}(p_{1, 1}) \\
& = \frac{1}{(g - 3)!}p_{2, 0}^{g - 3}p_{3, 1} - \frac{1}{(g - 2)!}p_{2, 0}^{g - 2}p_{1, 1}. \nonumber
\end{align*}

Now that $[C]_{(1)}$ is expressed in terms of $\{p_{i, j}\}$, we have
\begin{equation*}
[C]_{(1)} \cdot p_{3, 1} = \frac{(-1)^{g + 1}}{(g - 3)!}p_{2,0}^{g - 3}p_{3, 1}^2 - \frac{(-1)^{g + 1}}{(g - 2)!}p_{2, 0}^{g - 2}p_{1, 1}p_{3, 1}.
\end{equation*}
Apply $\mathcal{F}$ one more time to get the class we want:
\begin{align*}
-p_{3, 1} * [C]_{(1)} & = \mathcal{F}\bigg(\frac{(-1)^{g + 1}}{(g - 3)!}p_{2,0}^{g - 3}p_{3, 1}^2 - \frac{(-1)^{g + 1}}{(g - 2)!}p_{2, 0}^{g - 2}p_{1, 1}p_{3, 1}\bigg) \\
& = \exp(e)\exp(-f)\exp(e)\bigg(\frac{(-1)^{g + 1}}{(g - 3)!}e^{g - 3}(p_{3, 1}^2) - \frac{(-1)^{g + 1}}{(g - 2)!}e^{g - 2}(p_{1, 1}p_{3, 1})\bigg) \nonumber \\
& = \frac{(-1)^g}{g!}f^g\bigg(\frac{(-1)^{g + 1}}{(g - 3)!}e^{g - 3}(p_{3, 1}^2) - \frac{(-1)^{g + 1}}{(g - 2)!}e^{g - 2}(p_{1, 1}p_{3, 1})\bigg) \nonumber \\
& = -\frac{1}{g!(g - 3)!}f^ge^{g - 3}(p_{3, 1}^2) + \frac{1}{g!(g - 2)!}f^ge^{g - 2}(p_{1, 1}p_{3, 1}). \nonumber
\end{align*}
Expressions such as $f^ge^{g - 3}(p_{3, 1}^2)$ and $f^ge^{g - 2}(p_{1, 1}p_{3, 1})$ can be computed using a combinatorial formula for $\mathfrak{sl}_2$-representations. Here we state this formula as a lemma, since we shall use it again later. $\phantom{\qedhere}$
\end{proof}

\subsection{\it Lemma} \label{combi} --- {\it Consider a representation $\mathfrak{sl}_2 \to {\rm End}_{\mathbb{Q}}(V)$. Let $\alpha \in V$ such that $h(\alpha) = \mu \cdot \alpha$. Then for all $r, s \geq 0$ we have}
\begin{equation}
f^se^r(\alpha) = \sum_{t = 0}^{\min(r, s)}(-1)^t \frac{s!}{(s - t)!} \frac{r!}{(r - t)!} \binom{\mu + r - s + t - 1}{t} e^{r - t}f^{s - t}(\alpha).
\end{equation}

\medskip
Note that the binomial coefficient should be taken in the generalized sense. We refer to \cite{Moo09}, Lemma~2.4 for the proof (where the operator $\mathcal{D}$ in {\it loc.\,cit.}~is $-f$).

\begin{proof}[Proof of the Theorem {\rm (continued)}]
It follows from the lemma above that
\begin{equation*}
f^ge^{g - 3}(p_{3, 1}^2) = \frac{g!(g - 3)!}{3!}f^3(p_{3, 1}^2), \quad f^ge^{g - 2}(p_{1, 1}p_{3, 1}) = \frac{g!(g - 2)!}{2!}f^2(p_{1, 1}p_{3, 1}),
\end{equation*}
so that
\begin{equation*}
-p_{3, 1} * [C]_{(1)} = -\frac{1}{6}f^3(p_{3, 1}^2) + \frac{1}{2}f^2(p_{1, 1}p_{3, 1}).
\end{equation*}
By Theorem~\ref{structure} we know that $f = \mathcal{D}$ on polynomials in $\{p_{i, j}\}$ and $\psi$, so it remains to calculate
\begin{align*}
\mathcal{D}^3(p_{3, 1}^2) & = \mathcal{D}^2(\psi p_{0, 2}^2 - 6p_{4, 2} + 2p_{1, 1}p_{3, 1}) \\
& = \mathcal{D}\big((2g - 2)\psi p_{0, 2} - 6p_{2, 2} + 2(g \psi p_{0, 2} - p_{2, 2} + p_{1, 1}^2)\big) \nonumber \\
& = \mathcal{D}\big(2(2g - 1)\psi p_{0, 2} - 8p_{2, 2} + 2p_{1, 1}^2\big) \nonumber \\
& = 2g(2g - 1)\psi - 8p_{0, 2} + 2(g^2\psi - p_{0, 2}) \nonumber \\
& = 2g(3g - 1)\psi - 10p_{0, 2}, \nonumber
\end{align*}
and
\begin{align*}
\mathcal{D}^2(p_{1, 1}p_{3, 1}) & = \mathcal{D}(g\psi p_{2, 0} - p_{2, 2} + p_{1, 1}^2) \\
& = g^2\psi - p_{0, 2} + (g^2\psi - p_{0, 2}) \nonumber \\
& = 2g^2\psi - 2p_{0, 2}. \nonumber
\end{align*}
Altogether we have
\begin{equation*}
-p_{3, 1} * [C]_{(1)} = \frac{g}{3}\psi + \frac{2}{3}p_{0, 2} = \frac{g}{3}\psi + \frac{2}{3}\pi^*\bigg(\frac{1}{4}\kappa_1 + g\psi\bigg) = \pi^*\bigg(\frac{1}{6}\kappa_1 + g\psi\bigg),
\end{equation*}
where the middle equality follows from (\ref{basechange}). Combining with the starting point (\ref{goal}), we find
\begin{equation*}
\pi^*\pi_*\big([C]_{(1)} \cdot p_{3, 1}\big) = \pi^*\bigg(\frac{1}{6}\kappa_1 + g\psi\bigg).
\end{equation*}
Therefore $\pi_*\big([C]_{(1)} \cdot p_{3, 1}\big) = \kappa_1/6 + g\psi$ by the isomorphisms (\ref{tautisom}).
\end{proof}

\medskip
\section{\bf Relations between tautological classes: Faber's conjectures}
\medskip

\noindent Now that we know the generators ($\{p_{i, j}\}$ and $\psi$) of $\mathcal{T}(J)$, the question becomes `what are the relations between these classes?'. In this section we will see how the $\mathfrak{sl}_2$-action introduced earlier, together with its explicit form on $\mathcal{T}(J)$ ({\it cf.}~Theorem~\ref{structure}), produces relations between the generators. The idea goes back to Polishchuk in his paper \cite{Pol05}.

Thanks to Corollary~\ref{subring}, from now on we identify $\mathcal{R}(S)$ with $\oplus_{i = 0}^{d}\mathcal{T}_{(0, 2i)}(J)$ via the map $\pi^*$. Then by restricting everything to $\mathcal{R}(S)$, we obtain relations between the generators of $\mathcal{R}(S)$. This suggests a new approach to study the structure of $\mathcal{R}(\mathcal{M}_{g, 1})$, which is the subject of Faber's conjectures. We shall discuss various theoretical aspects and numerical evidence regarding these conjectures.

We start by stating Faber's conjectures properly in the context of $\mathcal{M}_{g, 1}$. See \cite{Fab99} for the original conjectures for $\mathcal{M}_g$ (moduli space of smooth curves of genus $g$).

\subsection{\it Conjectures {\rm (Faber, for \texorpdfstring{$\mathcal{M}_{g, 1}$}{Mg,1})}} \label{faber} --- (i) {\it The tautological ring $\mathcal{R}(\mathcal{M}_{g, 1})$ is Gorenstein with socle in codimension $g - 1$. This means that $\mathcal{R}^i(\mathcal{M}_{g, 1}) = 0$ for $i > g - 1$, that $\mathcal{R}^{g - 1}(\mathcal{M}_{g, 1}) \simeq \mathbb{Q}$, and that the natural paring}
\begin{equation} \label{pairing}
\mathcal{R}^i(\mathcal{M}_{g, 1}) \times \mathcal{R}^{g - 1 - i}(\mathcal{M}_{g, 1}) \xrightarrow{\textrm{`}\cdot\textrm{'}} \mathcal{R}^{g - 1}(\mathcal{M}_{g, 1}) \simeq \mathbb{Q}
\end{equation}
{\it is perfect for all $0 \leq i \leq g - 1$.}

\smallskip
\noindent\makebox[\leftmargin][r]{(ii) }{\it The ring  $\mathcal{R}(\mathcal{M}_{g, 1})$ is generated by $\kappa_1, \ldots, \kappa_{\lfloor g/3 \rfloor}$ and $\psi$. There are no relations between these classes in codimension $\leq \lfloor g/3 \rfloor$.}

\subsection{\it Remarks} \label{faberrem} --- (i) There are isomorphisms of stacks $\mathcal{M}_{g, 1} \simeq \mathcal{C}_g \simeq \mathcal{M}_{g, 1}^{\textrm{rt}}$, where $\mathcal{C}_g$ is the universal curve over $\mathcal{M}_g$, and $\mathcal{M}_{g, 1}^{\textrm{rt}}$ is the moduli space of stable $1$-pointed curves of genus $g$ with rational tails. Hence Faber's conjectures for these three moduli spaces are equivalent.

\smallskip
\noindent\makebox[\leftmargin][r]{(ii) }Looijenga \cite{Loo95} proved that $\mathcal{R}^i(\mathcal{M}_{g, 1})$ vanishes for $i > g - 1$, and that $\mathcal{R}^{g - 1}(\mathcal{M}_{g, 1})$ is 
at most $1$-dimensional. Later Faber himself \cite{Fab95} proved that $\mathcal{R}^{g - 1}(\mathcal{M}_{g, 1})$ is indeed $1$-dimensional. This means in the first conjecture, only the `perfect paring' statement remains open.

\smallskip
\noindent\makebox[\leftmargin][r]{(iii) }The generation statement in the second conjecture was first proven by Ionel \cite{Ion05}. (It follows from her result that the tautological ring of $\mathcal{M}_g$ is generated by $\kappa_1, \ldots, \kappa_{\lfloor g/3 \rfloor}$.) By (\ref{basechange}), saying that $\mathcal{R}(\mathcal{M}_{g, 1})$ is generated by $\kappa_1, \ldots, \kappa_{\lfloor g/3 \rfloor}$ and $\psi$ is equivalent to saying that it is generated by $p_{0, 2}, \ldots, p_{0, 2\lfloor g/3 \rfloor}$ and $\psi$ ({\it cf.}~Corollary~\ref{subring}).

\smallskip
\noindent\makebox[\leftmargin][r]{(iv) }There should be a third conjecture similar to the one for $\mathcal{M}_g$, which gives explicit formulae for the proportionalities in codimension $g - 1$. This could be obtained in an {\it ad hoc\/} manner by pushing forward to $\mathcal{M}_g$ ({\it cf.}~Section~\ref{Mg}). See Section~\ref{socle} for some other discussion.

\smallskip
\noindent\makebox[\leftmargin][r]{(v) }Faber also considered a stronger version of the first conjecture. In the case of $\mathcal{M}_{g, 1}$, it says that $\mathcal{R}(\mathcal{M}_{g, 1})$ `behaves like' the algebraic cohomology ring of a smooth projective variety of dimension $g - 1$, {\it i.e.}~it satisfies the hard Lefschetz and Hodge positivity properties. From this perspective, the generators $\{p_{0, 2i}\}$ and $\psi$ have some advantage over $\{\kappa_i\}$ and $\psi$: both $p_{0, 2}$ and $\psi$ are candidates of the conjectural `ample' divisor class for hard Lefschetz and Hodge positivity, while on the other hand we have $\kappa_1^{g - 1} = 0$.

\subsection{} \label{relation} We explain how the $\mathfrak{sl}_2$-action produces relations between the generators of $\mathcal{T}(J)$ (resp.~$\mathcal{R}(S)$). First, Theorem~\ref{structure} shows that the space of polynomial relations between the classes $\{p_{i, j}\}$ and $\psi$ is stable under the action of $\mathcal{D}$. In other words, for all polynomials $P$ in $\{p_{i, j}\}$ and $\psi$, we have that $P\big(\{p_{i, j}\}, \psi\big) = 0$ implies $\mathcal{D}\Big(P\big(\{p_{i, j}\}, \psi\big)\Big) = 0$. Then consider monomials
\begin{equation*}
\alpha = \psi^s p_{i_1, j_1}^{r_1} p_{i_2, j_2}^{r_2} \cdots p_{i_m, j_m}^{r_m}\, \textrm{ with } I := r_1 i_1 + r_2 i_2 + \cdots + r_m i_m > 2g.
\end{equation*}
By definition $\alpha \in {\rm CH}\big(R^I(J/S)\big)$. But since $I > 2g$, we know from decomposition~(\ref{motive})
that $R^I(J/S) = 0$. In terms of the Dutch house, the class $\alpha$ is simply outside the house. It follows that we have relations
\begin{equation*}
\alpha = 0, \quad \mathcal{D}(\alpha) = 0, \quad \mathcal{D}^2(\alpha) = 0, \quad \ldots
\end{equation*}

This argument leads to the following formal definition.

\subsection{} \label{abs} Let $i, j$ run through all integers such that $i \leq j + 2$ and that $i + j$ is even. Consider the ring
\begin{equation*}
\mathcal{A} := \mathbb{Q}\big[\{x_{i, j}\}, y\big]\left/\big<x_{0, 0} - g, \{x_{i, j}\}_{i < 0}, \{x_{i, j}\}_{j < 0}, \{x_{i, j}\}_{j > 2g - 2}\big>.\right.
\end{equation*}
In other words, the ring $\mathcal{A}$ is a polynomial ring in variables $\{x_{i, j}\}$ and $y$, with the convention that $x_{0, 0} = g$ and $x_{i, j} = 0$ for $i < 0 $ or $j < 0$ or $j > 2g - 2$ ({\it cf.}~Section~\ref{pij}). We introduce a bigrading $\mathcal{A} = \oplus_{i, j}\mathcal{A}_{(i, j)}$ by the requirements that $x_{i, j} \in \mathcal{A}_{(i, j)}$ and $y \in \mathcal{A}_{(0, 2)}$. Define operators $E, F$ and $H$ on $\mathcal{A}$ by
\begin{align*}
E(\alpha) & := x_{2, 0} \cdot \alpha, \\
F(\alpha) & :=  \frac{1}{2}\sum_{\substack{i, j \\ k, l}}\Bigg(y x_{i - 1, j - 1} x_{k - 1, l - 1} - \binom{i + k - 2}{i - 1} x_{i + k - 2, j + l}\Bigg) \partial x_{i, j} \partial x_{k, l}(\alpha) \\
& \qquad + \sum_{i, j}x_{i - 2, j}\partial x_{i, j}(\alpha), \nonumber \\
H(\alpha) & := (i - g) \cdot \alpha, \textrm{ for } \alpha \in \mathcal{A}_{(i, j)}.
\end{align*}
It is not difficult to verify that the operators above generate a representation $\mathfrak{sl}_2 \to \textrm{End}_{\mathbb{Q}}(\mathcal{A})$. Theorem~\ref{structure} can then be reformulated by saying that we have a surjective morphism of $\mathfrak{sl}_2$-representations $\mathcal{A} \to \mathcal{T}(J)$, which maps $x_{i, j}$ to $p_{i, j}$ and $y$ to $\psi$.

Denote by ${\rm Mon}_{(i, j)}$ the set of monomials in $\{x_{i, j}\}$ (excluding $x_{0, 0}$ and all $x_{i, j}$ that vanish in $\mathcal{A}$) and $y$ that belong to $\mathcal{A}_{i, j}$. Note that we set ${\rm Mon}_{(0, 0)} = \{1\}$ as an exception. Then consider the quotient ring
\begin{equation}
\tilde{\mathcal{T}} := \mathcal{A}\left/\Big<\big\{F^\nu({\rm Mon}_{(i, j)})\big\}_{i > 2g, \nu \geq 0}\Big>.\right.
\end{equation}
\noindent The ring $\tilde{\mathcal{T}}$ inherits from $\mathcal{A}$ a bigrading $\tilde{\mathcal{T}} = \oplus_{i, j}\tilde{\mathcal{T}}_{(i, j)}$. The operators $E, F$ and $H$ induce operators on $\tilde{\mathcal{T}}$, which we denote by $e, f$ and $h$. Again we obtain a representation $\mathfrak{sl}_2 \to \textrm{End}_{\mathbb{Q}}(\tilde{\mathcal{T}})$. Moreover, since $e^{g + 1} = f^{g + 1} = 0$, we can formally define the Fourier transform $\mathcal{F}$ on $\tilde{\mathcal{T}}$ by
\begin{equation*}
\mathcal{F} := \exp(e)\exp(-f)\exp(e).
\end{equation*}

\begin{figure}
\centering
\includegraphics[height=.43\textheight]{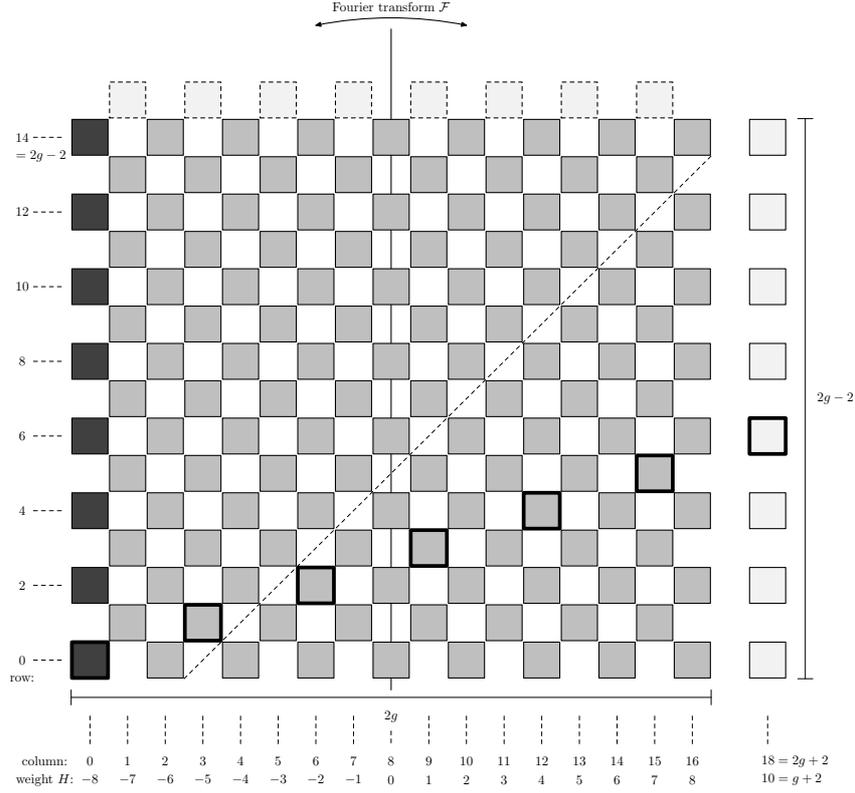}
\caption{Producing relations ($g = 8$).}
\label{relations}
\end{figure}

We define also the subring
\begin{equation}
\tilde{\mathcal{R}} = \oplus_i \tilde{\mathcal{T}}_{(0, 2i)},
\end{equation}
with the grading $\tilde{\mathcal{R}} = \oplus_i \tilde{\mathcal{R}}^i$ such that $\tilde{\mathcal{R}}^i := \tilde{\mathcal{T}}_{(0, 2i)}$. Then we have
\begin{align*}
\tilde{\mathcal{R}}^i & = \mathcal{A}_{(0, 2i)}\left/\Big<\big\{F^I({\rm Mon}_{(2I, 2i)})\big\}_{I > g}\Big>\right.\\
& = \mathcal{A}_{(0, 2i)}\left/\big<F^{g + 1}({\rm Mon}_{(2g + 2, 2i)})\big>.\right. \nonumber
\end{align*}
Figure~\ref{relations} illustrates the construction of $\tilde{\mathcal{R}}$: take monomials on the $(2g + 2)$-th column of the Dutch house (white blocks), then apply $g + 1$ times the operator $F$ to obtain relations between the generators (black blocks).

To summarize this formal approach, we have the following proposition.

\subsection{\it Proposition} \label{surjections} --- {\it For all $S$ as in diagram~{\rm (\ref{setting}) (}including $S = \mathcal{M}_{g, 1}${\rm )}, we have surjective maps
\begin{equation}
\Phi \colon \tilde{\mathcal{T}} \twoheadrightarrow \mathcal{T}(J), \quad \Phi|_{\tilde{\mathcal{R}}} \colon \tilde{\mathcal{R}} \twoheadrightarrow \mathcal{R}(S),
\end{equation}
which map $x_{i, j}$ to $p_{i, j}$ and $y$ to $\psi$.}

\medskip
From now on, we will concentrate on the structure of $\tilde{\mathcal{R}}$. We start with a lemma, which shows we can already eliminate certain monomials that produce trivial relations.

\subsection{\it Lemma} --- {\it For all $\alpha \in {\rm Mon}_{(2g + 2, 2i)}$ of the form $\alpha = x_{2, 0} \cdot \beta$, we have $F^{g + 1}(\alpha) = 0$.}

\begin{proof}
This follows directly from Lemma~\ref{combi}. In fact, we have $F^{g + 1}(\alpha) = F^{g + 1}(x_{2, 0} \cdot \beta) = F^{g + 1}E(\beta)$. Then apply Lemma~\ref{combi} with $\mu = g$, $r = 1$ and $s = g + 1$, and we find
\begin{equation*}
F^{g + 1}E(\beta) = EF^{g + 1}(\beta) - (g + 1)(g + 1 - g - 1 + 1 - 1)F^g(\beta) = EF^{g + 1}(\beta).
\end{equation*}
On the other hand $F^{g + 1}(\beta) \in \mathcal{A}_{(-2, 2i)} = 0$, which implies $F^{g + 1}(\alpha) = 0$.
\end{proof}

As a result, if we write ${\rm mon}_{(2g + 2, 2i)} \subset {\rm Mon}_{(2g + 2, 2i)}$ for the subset of monomials without $x_{2, 0}$ as a factor, then we have
\begin{equation} \label{defR}
\tilde{\mathcal{R}}^i = \mathcal{A}_{(0, 2i)}\left/\big<F^{g + 1}({\rm mon}_{(2g + 2, 2i)})\big>.\right. 
\end{equation}
The bold blocks in Figure~\ref{relations} describe the lower bound of $i$ such that ${\rm mon}_{(2g + 2, 2i)}$ is non-empty. In numerical terms, we have $x_{3, 1}^{2i} \in \mathcal{A}_{(6i, 2i)}$, and then $6i > 2g$ implies $i > g/3$. It follows that ${\rm mon}_{(2g + 2, 2i)} = \emptyset$ for all $i \leq \lfloor g/3 \rfloor$.

\subsection{\it Theorem} \label{generation} --- {\it The ring $\tilde{\mathcal{R}}$ is generated by $x_{0, 2}, \cdots, x_{0, 2 \lfloor g/3 \rfloor}$ and $y$. There are no relations between these elements in $\tilde{\mathcal{R}}^i$ for $i \leq \lfloor g/3 \rfloor$.}

\medskip
In particular, this gives a new proof of the generation statement of Conjecture~\ref{faber}~(ii) ({\it cf.}~Remark~\ref{faberrem}~(iii)). Further, it proves Conjecture~\ref{faber}~(ii) completely if the map $\Phi|_{\tilde{\mathcal{R}}} \colon \tilde{\mathcal{R}} \twoheadrightarrow \mathcal{R}(\mathcal{M}_{g, 1})$ in Proposition~\ref{surjections} is also injective, {\it i.e.}~an isomorphism.

\begin{proof}
The second statement is immediate after (\ref{defR}) and the fact that ${\rm mon}_{(2g + 2, 2i)} = \emptyset$ for all $i \leq \lfloor g/3 \rfloor$. For the first statement, the goal is to relate all $x_{0, 2i}$ with $i > g/3$ to the elements $x_{0, 2}, \cdots, x_{0, 2 \lfloor g/3 \rfloor}$ and $y$, and the idea is to use specific monomials to get these relations.

We proceed by induction. Suppose all $\{x_{0, 2j}\}_{g/3 < j < i}$ can be expressed in terms of $x_{0, 2}, \cdots, x_{0, 2 \lfloor g/3 \rfloor}$ and $y$. Then consider the monomial $x_{3, 1}^{2i} \in \mathcal{A}_{(6i, 2i)}$. Apply $3i$-times the operator $F$ and we get $F^{3i}(x_{3, 1}^{2i}) \in \mathcal{A}_{(0, 2i)}$, which vanishes in $\tilde{\mathcal{R}}$. By the explicit formula of $F$, we have
\begin{equation*}
F^{3i}(x_{3, 1}^{2i}) = c \cdot x_{0, 2i} + \alpha,
\end{equation*}
where $\alpha$ is a polynomial in $\{x_{0, 2j}\}_{j < i}$ and $y$. It only remains to prove that the coefficient $c$ is non-zero.

The observation is the following: when we apply the operator $F$, the minus sign occurs every time two factors ($x_{i, j}, x_{k, l}$) are merged into one ($x_{i + k - 2, j + l}$). Then if we start from 
$x_{3, 1}^{2i}$ and arrive at $x_{0, 2i}$, no matter how we proceed we have to do the merging $(2i - 1)$-times. This means that all non-zero summands of $c$ are of the form $(-1)^{2g - 1}$ times a positive number, hence negative. Therefore the sum $c$ is negative as well.
\end{proof}

\subsection{} \label{socle} Consider the component $\tilde{\mathcal{R}}^{g - 1}$, where lives the expected socle. It would be nice if we could prove $\tilde{\mathcal{R}}^{g - 1} \simeq \mathbb{Q}$ using our relations. We already knew that $\tilde{\mathcal{R}}^{g - 1}$ is spanned by the image of ${\rm Mon}_{(0, 2g - 2)}$. The goal is then to relate all elements of ${\rm Mon}_{(0, 2g - 2)}$ to a single element, namely $x_{0, 2g - 2}$. So far we have been able to associate a relation to every element of ${\rm Mon}_{(0, 2g - 2)}$ other than $x_{0, 2g - 2}$. By defining a partial order on  ${\rm Mon}_{(0, 2g - 2)}$, we can show that the relation matrix is block triangular. It remains to prove that all diagonal blocks are of maximal rank, which is combinatorially difficult.

The construction goes as follows: first remark that there is a correspondence between ${\rm Mon}_{(0, 2g - 2)}$ and the set of partitions of all integers $k$ such that $0 \leq k \leq g - 1$. In fact, to every such partition $\lambda = (i_1, \ldots, i_m)$ with $i_1 + \cdots + i_m = k$, we associate the element $y^{g - 1 - k} x_{0, 2i_1} \cdots x_{0, 2i_m}$ in ${\rm Mon}_{(0, 2g - 2)}$. In particular, we have $\#{\rm Mon}_{(0, 2g - 2)} = p(0) + \cdots + p(g - 1)$, where $p(\textrm{-})$ stands for the partition function. Then observe that to every partition $\lambda$ above, we may associate a partition $\lambda'$ of $2g - 2$
\begin{equation*}
\lambda = (i_1, \ldots, i_m) \mapsto \lambda' = (2i_1 + g - 1 - k, 2i_2, \ldots, 2i_m, \underbrace{1, \ldots, 1}_{g - 1 - k}).
\end{equation*}
The partial order on ${\rm Mon}_{(0, 2g - 2)}$ is defined as the first part of $\lambda'$, {\it i.e.}~$2i_1 + g - 1 - k$. Using the elements $\{x_{i + 2, i}\}$ located to the right of the diagonal $i = j$ in the Dutch house, we define the monomial
\begin{equation*}
M_{\lambda'} = x_{2i_1 + g - 1 - k + 2, 2i_1 + g - 1 - k} \cdot x_{2i_2 + 2, 2i_2} \cdots x_{2i_m + 2, 2i_m} \cdot x_{3, 1}^{g - 1 - k},
\end{equation*}
which belongs to ${\rm Mon}_{(4g - 4 + 2m - 2k, 2g - 2)}$. We see that for all $\lambda$ except $\lambda = (g - 1)$, we have $4g - 4 + 2m - 2k > 2g$, so that $M_{\lambda'} = 0$ in $\tilde{\mathcal{R}}$. Then to every element of ${\rm Mon}_{(0, 2g - 2)}$ other than $x_{0, 2g - 2}$, we may associate a relation $F^{2g - 2 + m - k}(M_{\lambda'}) = 0$ in $\tilde{\mathcal{R}}^{g - 1}$. By the explicit formula of $F$, we can prove that $F^{2g - 2 + m - k}(M_{\lambda'})$ only contains terms of order greater than or equal to the first part of $\lambda'$.

For $\tilde{\mathcal{R}}^{i}$ with $i > g - 1$, we are in similar situation. In this case there are no more new generators $x_{0, 2i}$, and every element in ${\rm Mon}_{(0, 2i)}$ corresponds to a relation. The question is still to show that the relation matrix is of maximal rank. Here we should say that although the combinatorics seem too difficult for the moment, we do get a feeling why the socle should lie in codimension $g - 1$, and why beyond codimension $g - 1$ everything should vanish.

\subsection{\it Numerical evidence} --- Li~Ma made a C++ program that computes the ring $\tilde{\mathcal{R}}$ for a given genus $g$. The program calculates all relations and outputs the dimension of each component $\tilde{\mathcal{R}}^i$. We refer to the appendix for more computational  details.

On the other hand, based on an algorithm developed by Liu and Xu \cite{LX12}, Olof~Bergvall \cite{Ber11} computed the proportionalities of the classes in $\mathcal{R}^{g - 1}(\mathcal{M}_{g, 1}) \simeq \mathbb{Q}$ for many values of $g$. As a result, it gives the dimensions of the Gorenstein quotient $\mathcal{G}(\mathcal{M}_{g, 1})$ of $\mathcal{R}(\mathcal{M}_{g, 1})$. Here the ring $\mathcal{G}(\mathcal{M}_{g, 1})$ is the quotient of $\mathcal{R}(\mathcal{M}_{g, 1})$ by the homogeneous ideal generated by all classes of pure codimension that have zero pairing with all classes of the opposite codimension ({\it cf.}~\cite{Fab11}). Note that this procedure of obtaining the dimensions of $\mathcal{G}(\mathcal{M}_{g, 1})$ is formal and does not involve computation of actual relations in $\mathcal{R}(\mathcal{M}_{g, 1})$.

There are surjective maps $\tilde{\mathcal{R}} \twoheadrightarrow \mathcal{R}(\mathcal{M}_{g, 1}) \twoheadrightarrow \mathcal{G}(\mathcal{M}_{g, 1})$. Our computation shows that for $g \leq 19$, the dimensions of $\tilde{\mathcal{R}}$ and $\mathcal{G}(\mathcal{M}_{g, 1})$ are the same, which implies isomorphisms $\tilde{\mathcal{R}} \simeq \mathcal{R}(\mathcal{M}_{g, 1}) \simeq \mathcal{G}(\mathcal{M}_{g, 1})$. In particular, we can confirm the following. (For $g \leq 9$ this has been obtained independently by Bergvall \cite{Ber11}.)

\subsection{\it Theorem} \label{g19} --- {\it Conjectures~\ref{faber} are true for $g \leq 19$.}

\medskip
However, the computer result is negative for $g = 20$ and some greater values of $g$. There the dimensions of $\tilde{\mathcal{R}}$ are simply not symmetric. Again by comparing with the dimensions of $\mathcal{G}(\mathcal{M}_{g, 1})$, we know exactly how many relations are missing. The numbers are listed in Table~\ref{tab} below.

\subsection{} \label{sequence} The computation above also involves an observation similar to Faber and Zagier's approach. Consider the dimension of $\tilde{\mathcal{R}}^i$ for $g/3 < i \leq \big\lfloor (g - 1)/2 \big\rfloor$, which we denote by $d(g, i)$. Note that $i \leq \big\lfloor (g - 1)/2 \big\rfloor$ corresponds to the lower half of the conjectural Gorenstein ring. We know that $\tilde{\mathcal{R}}^i$ is spanned by the image of ${\rm Mon}_{(0, 2i)}$, and by an argument in Section~\ref{socle} we have $\#{\rm Mon}_{(0, 2i)} = p(0) + \cdots + p(i) =: \phi(i)$. Then consider the difference $\phi(i) - d(g, i)$, which reflects the minimal number of relations needed to obtain $\tilde{\mathcal{R}}^i$ from $\mathcal{A}_{(0, 2i)}$ ({\it cf.}~(\ref{defR})). We have observed an interesting phenomenon: that $\phi(i) - d(g, i)$ seems to be a function of $3i - g - 1$. More precisely, if we set
\begin{equation*}
b(3i - g - 1) := \phi(i) - d(g, i),
\end{equation*}
then based on our data for $g \leq 28$, we obtain the following values of the function $b$.
\smallskip
\begin{center}
{\renewcommand{\arraystretch}{1.15}
\begin{tabular}{c||c|c|c|c|c|c|c|c|c|c|c|c}
\hline
$n$ & $0$ & $1$ & $2$ & $3$ & $4$ & $5$ & $6$ & $7$ & $8$ & $9$ & $10$ & $11$ \\
\hline
$b(n)$ & $1$ & $2$ & $3$ & $6$ & $10$ & $14$ & $22$ & $33$ & $45$ & $64$ & $90$ & $119$ \\
\hline
\end{tabular}}
\end{center}

\smallskip
In \cite{Ber11}, Bergvall did the same calculation with the dimensions of the Gorenstein quotient $\mathcal{G}(\mathcal{M}_{g, 1})$. There he obtained the same values as the function $b(n)$ for $n \leq 9$. For $n = 10$, however, he observed some unexpected fact: the component $\mathcal{G}^{12}(\mathcal{M}_{25, 1})$ yields the number $91$ instead of $90$, while further components $\mathcal{G}^{13}(\mathcal{M}_{28, 1})$ and $\mathcal{G}^{14}(\mathcal{M}_{31, 1})$ give the number $90$. Similar things happen to $n = 11$, where the component $\mathcal{G}^{13}(\mathcal{M}_{27, 1})$ yields the number $120$ instead of $119$, while further components seem to give the number $119$. Note that the `defect' occurs only in the middle codimension, {\it i.e.}~$(g - 1)/2$ for $g$ odd.

Another aspect is to guess what this function $b(n)$ could be. Recall Faber and Zagier's function $a(n)$ for $\mathcal{M}_g$ ({\it cf.}~\cite{Fab99} and \cite{LX12}).
\smallskip
\begin{center}
{\renewcommand{\arraystretch}{1.15}
\begin{tabular}{c||c|c|c|c|c|c|c|c|c|c|c|c}
\hline
$n$ & $0$ & $1$ & $2$ & $3$ & $4$ & $5$ & $6$ & $7$ & $8$ & $9$ & $10$ & $11$ \\
\hline
$a(n)$ & $1$ & $1$ & $2$ & $3$ & $5$ & $6$ & $10$ & $13$ & $18$ & $24$ & $33$ & $41$ \\
\hline
\end{tabular}}
\end{center}
\smallskip
It has been suggested by Bergvall and Faber that the first values of $b(n)$ satisfy
\begin{equation} \label{BF}
b(n) = \sum_{\substack{i = 0 \\ i \not\equiv 2 \,({\rm mod} \,3)}}^n a(n - i) = a(n) + a(n - 1) + a(n - 3) + a(n - 4) + \cdots.
\end{equation}
It looks mysterious that Formula~(\ref{BF}) is compatible with the dimensions of $\tilde{\mathcal{R}}$, but conflicts with the Gorenstein conjecture.

\subsection{\it Pushing forward to \texorpdfstring{$\mathcal{M}_g$}{Mg}} \label{Mg} --- The tautological ring $\mathcal{R}(\mathcal{M}_g)$ is defined to be the $\mathbb{Q}$-subalgebra of ${\rm CH}(\mathcal{M}_g)$ generated by the classes $\{\kappa_i\}$. Faber's conjectures in this case predict that $\mathcal{R}(\mathcal{M}_g)$ is Gorenstein with socle in codimension $g - 2$, that it is generated by $\kappa_1, \ldots, \kappa_{\lfloor g/3 \rfloor}$, and that there are no relations between these classes in codimension $\leq \lfloor g/3 \rfloor$. See \cite{Fab99} for more details.

Let $q \colon \mathcal{M}_{g, 1} \to \mathcal{M}_g$ be the map that forgets the marked point. Then we have $q_*\big(\mathcal{R}(\mathcal{M}_{g, 1})\big) = \mathcal{R}(\mathcal{M}_g)$. In fact, consider the following commutative diagram
\begin{equation*}
\begin{tikzcd}[row sep=large, column sep=large]
\mathcal{C}_{g, 1} \arrow{r}{\tilde{q}} \arrow{d}[swap]{\tilde{p}} & \mathcal{C}_g \arrow{d}{p} \\
\mathcal{M}_{g, 1} \arrow{r}[swap]{q} \arrow[bend left=50]{u}{x_0} & \mathcal{M}_g
\end{tikzcd}
\end{equation*}
where $\mathcal{C}_g$ (resp.~$\mathcal{C}_{g, 1}$) is the universal curve over $\mathcal{M}_g$ (resp.~$\mathcal{M}_{g, 1}$), and $x_0$ is the section of $\tilde{p}$ that gives the marked point. The classes $\{\kappa_i\}$ (resp.~$K$) are defined on $\mathcal{M}_g$ (resp.~$\mathcal{C}_g$), and we use the same notation for their pullback to $\mathcal{M}_{g, 1}$ (resp.~$\mathcal{C}_{g, 1}$). Then for all $\psi^s \kappa_{i_1}^{r_1} \cdots \kappa_{i_m}^{r_m} \in \mathcal{R}(\mathcal{M}_{g, 1})$, we have
\begin{align} \label{pushforward}
q_*(\psi^s \kappa_{i_1}^{r_1} \cdots \kappa_{i_m}^{r_m}) & = q_*\big(\psi^s \cdot q^*(\kappa_{i_1}^{r_1} \cdots \kappa_{i_m}^{r_m})\big) \\
& = q_*(\psi^s) \cdot \kappa_{i_1}^{r_1} \cdots \kappa_{i_m}^{r_m} \nonumber \\
& = q_*x_0^*(K^s) \cdot \kappa_{i_1}^{r_1} \cdots \kappa_{i_m}^{r_m} \nonumber \\
& = q_*x_0^*\,\tilde{q}^*(K^s) \cdot \kappa_{i_1}^{r_1} \cdots \kappa_{i_m}^{r_m} \nonumber \\
& = p_*(\tilde{q} \circ x_0)_*(\tilde{q} \circ x_0)^*(K^s) \cdot \kappa_{i_1}^{r_1} \cdots \kappa_{i_m}^{r_m} \nonumber \\
& = p_*(K^s) \cdot \kappa_{i_1}^{r_1} \cdots \kappa_{i_m}^{r_m} \tag{$(\tilde{q} \circ x_0)_*(\mathbbm{1}_{\mathcal{M}_{g, 1}}) = \mathbbm{1}_{\mathcal{C}_g}$} \\
& = \kappa_{s - 1} \kappa_{i_1}^{r_1} \cdots \kappa_{i_m}^{r_m}, \nonumber
\end{align}
with the convention that $\kappa_{-1} = 0$. Hence $q_*(\psi^s \kappa_{i_1}^{r_1} \cdots \kappa_{i_m}^{r_m}) \in \mathcal{R}(\mathcal{M}_g)$.

Formulae~(\ref{basechange}) and (\ref{pushforward}) allow us to push relations in $\mathcal{R}(\mathcal{M}_{g, 1})$ forward to $\mathcal{R}(\mathcal{M}_g)$. We used another computer program to do the work. Then for $g \leq 23$, we obtain a new proof of the following result of Faber-Zagier, Faber and Pandharipande-Pixton ({\it cf.}~\cite{Fab11}).

\subsection{\it Theorem {\rm (Faber-Zagier, Faber and Pandharipande-Pixton)}} \label{g23} --- {\it For $g \leq 23$, the ring $\mathcal{R}(\mathcal{M}_g)$ is Gorenstein with socle in codimension $g - 2$.}

\medskip
In this case, the rest of Faber's conjectures is also true by explicit calculation. Note that when $20 \leq g \leq 23$, the missing relations in $\mathcal{R}(\mathcal{M}_{g, 1})$ do not affect the Gorenstein property of $\mathcal{R}(\mathcal{M}_g)$. From $g = 24$ on, however, the computer result is again negative. Our computation for $g \leq 28$ suggests that we obtain exactly the same set of relations as the Faber-Zagier relations ({\it cf.}~\cite{Fab11} and \cite{PP11}). Notably in the crucial case of $g = 24$, we have not found the missing relation in codimension $12$. It is not yet known whether in theory we obtain the same relations. 

We summarize in Table~\ref{tab} the computer result for both $\mathcal{M}_{g, 1}$ and $\mathcal{M}_g$. Note that all codimensions have been calculated for $g \leq 24$. For $g \geq 25$, we only calculated a range near the middle codimension, presuming that the tautological rings behave well near the socle degree.

\begin{table}
\tabcolsep = 10pt
\begin{center}
{\renewcommand{\arraystretch}{1.15}
\begin{tabular}{c|r@{: }l|r@{: }l}
\hline
$g$ &  \multicolumn{2}{c|}{$\mathcal{M}_{g, 1}$} & \multicolumn{2}{c}{$\mathcal{M}_g$} \\
\hline
$\leq 19$ & \multicolumn{2}{c|}{OK} &  \multicolumn{2}{c}{OK} \\
20 & codim $10$ & $1$ missing &  \multicolumn{2}{c}{OK} \\
21 & codim $11$ & $1$ missing &  \multicolumn{2}{c}{OK} \\
22 & codim $11$ & $1$ missing &  \multicolumn{2}{c}{OK} \\
23 & codim $12$ & $3$ missing &  \multicolumn{2}{c}{OK} \\
24 & codim $13$ & $2$ missing & codim $12$ & $1$ missing \\
&  codim $12$ & $4$ missing \\
25 & codim $13$ & $5$ missing & codim $12$ & $1$ missing \\
& codim $12$ & $1$ missing \\
26 & codim $14$ & $6$ missing & codim $13$ & $1$ missing \\
& codim $13$ & $6$ missing \\
27 & codim $15$ & $3$ missing & codim $14$ & $1$ missing \\
& codim $14$ & $11$ missing & codim $13$ & $1$ missing \\
& codim $13$ & $1$ missing \\
28 & codim $15$ & $10$ missing & codim $14$ & $2$ missing \\
& codim $14$ & $10$ missing \\
\hline
\end{tabular}}
\end{center}
\medskip
\caption{Computer output for $g \leq 28$.}
\label{tab}
\end{table}

\subsection{\it Final speculation} --- Now that the ring $\tilde{\mathcal{R}}$ is in general not Gorenstein, there are two main possibilities.

\smallskip
\noindent\makebox[\leftmargin][r]{(i) }{\it The ring $\mathcal{R}(\mathcal{M}_{g, 1})$ is still Gorenstein.} In this case we need to find other ways that produce the missing relations. Earlier works ({\it cf.}~\cite{Moo09}, Section 2.13) suggest that our approach is closely related to Brill-Noether theory, while the relation is not yet clear. However, it seems that a new kind of geometry will be needed to prove Faber's conjectures.

\smallskip
\noindent\makebox[\leftmargin][r]{(ii) }{\it The ring $\mathcal{R}(\mathcal{M}_{g, 1})$ is isomorphic to $\tilde{\mathcal{R}}$.} There are several reasons to believe this. First, in terms of numbers our method produces huge quantities of relations ({\it cf.}~Appendix), much more than any of the previous methods. Yet the rank of the relation matrices stays the same. Second, it is believed that on a generic Jacobian variety, the $\mathfrak{sl}_2$-action is the only source of relations between tautological classes ({\it cf.}~\cite{Pol05}). The idea goes back to Beauville \cite{Bea04} and Ceresa \cite{Cer83}, and since then there have been various approaches to proving the non-triviality of certain tautological classes. We refer to Fakhruddin's paper \cite{Fak96} and some very recent work of Voisin (unpublished). It would be nice if we could find similar arguments for the universal Jacobian over $\mathcal{M}_{g, 1}$.

\smallskip
For the moment, we find it worthwhile to re-state the second case as a conjecture, although it contradicts Faber's conjectures for both $\mathcal{M}_{g, 1}$ and $\mathcal{M}_g$.

\subsection{\it Conjecture} \label{conj} --- {\it We have an isomorphism $\Phi|_{\tilde{\mathcal{R}}} \colon \tilde{\mathcal{R}} \xrightarrow{\sim} \mathcal{R}(\mathcal{M}_{g, 1})$.}

\medskip
\section*{\bf Appendix}
\begin{center}
{\Small LI MA}
\end{center}
\medskip

\noindent We briefly describe the computer program that calculates the dimensions of the $\mathbb{Q}$-vector spaces $\tilde{\mathcal{R}}^i$ introduced in Section~\ref{abs}.

\subsection*{\it {\rm I. }The algorithm} --- By (\ref{defR}), the vector space $\tilde{\mathcal{R}}^i$ is the quotient space of $\mathcal{A}_{(0, 2i)}$ by the subspace generated by the elements $F^{g + 1}(\alpha)$, for $\alpha$ running through the set ${\rm mon}_{(2g + 2, 2i)}$. Accordingly, the algorithm does the following.

\smallskip
\noindent\makebox[\leftmargin][r]{(i) }List all the elements of ${\rm Mon}_{(0, 2i)}$, which form a basis of the vector space $\mathcal{A}_{(0, 2i)}$ (call them `generators').

\smallskip
\noindent\makebox[\leftmargin][r]{(ii) }For every monomial in the set ${\rm mon}_{(2g + 2, 2i)}$, apply $(g + 1)$-times the operator $F$ to it; the resulting polynomial must be a linear combination of the generators (with coefficients in $\mathbb{Z}$), and we thus get a linear relation among the generators by letting the polynomial be zero.

\smallskip
\noindent\makebox[\leftmargin][r]{(iii) }Having done the second step for all the monomials in ${\rm mon}_{(2g + 2,2i)}$, we get a `relation matrix'; the dimension of $\tilde{\mathcal{R}}^i$ is then equal to the number of generators minus the rank of the relation matrix.

\smallskip
However, there are two major difficulties. First, our method produces too many relations, {\it i.e.}~the cardinality $\#{\rm mon}_{(2g + 2, 2i)}$ grows too fast with respect to $g$ and $i$. Therefore it is often not possible to calculate all of them, especially for $i$ close to $g - 1$. As an example, here is a list of the numbers of generators and relations for $g = 24$.

\medskip
{\renewcommand{\arraystretch}{1.15} \small
\begin{tabular}{c||c|c|c|c|c|c|c|c|c|c|c|c|c|c|c|c}
\hline
$i$ & $0$ & $1$ & $2$ & $3$ & $4$ & $5$ & $6$ & $7$ & $8$ & $9$ & $10$ & $11$ & $12$ & $13$ & $14$ & $15$ \\
\hline
$\#{\rm Mon}_{(0, 2i)}$ & $1$ & $2$ & $4$ & $7$ & $12$ & $19$ & $30$ & $45$ & $67$ & $97$ & $139$ & $195$ & $272$ & $373$ & $508$ & $684$ \\
\hline
$\#{\rm mon}_{(50, 2i)}$ & $0$ & $0$ & $0$ & $0$ & $0$ & $0$ & $0$ & $0$ & $0$ & $5$ & $49$ & $325$ & $1709$ & $7763$ & $31530$ & $117275$ \\
\hline 
\end{tabular}

\smallskip
\begin{tabular}{c||c|c|c|c|c|c|c|c}
\hline
$i$ & $16$ & $17$ & $18$ & $19$ & $20$ & $21$ & $22$ & $23$ \\
\hline
$\#{\rm Mon}_{(0, 2i)}$ & $915$ & $1212$ & $1597$ & $2087$ & $2714$ & $3506$ & $4508$ & $5763$ \\
\hline
$\#{\rm mon}_{(50, 2i)}$ & $404905$ & $1310010$ & $3995122$ & $11532380$ & $31602373$ & $82422889$ & $205123969$ & $488481821$ \\
\hline
\end{tabular}}

\medskip
In reality, since we know the dimensions of the Gorenstein quotient $\mathcal{G}(\mathcal{M}_{g, 1})$ from \cite{Ber11}, we stop the calculation as soon as the expected dimensions are reached. However, when there are missing relations ({\it cf.}~Table~\ref{tab}), we have to calculated all possible relations. Fortunately this only occurs near the middle codimension, where the number of relations is relatively small. For instance in genus $24$ and codimension $13$, we calculated all $7763$ relations, which gave a rank that was $2$ short of the expected one. When pushing forward to $\mathcal{M}_{24}$, there was again $1$ missing relation in codimension $12$.

The second difficulty comes from the fact that the coefficients of the relations are extremely large. One could use multiple precision to get accurate results, but neither the time complexity nor the space complexity would be tolerable. For similar reasons floating point numbers are not useful either.

Our choice is to take a large prime number $p$ and do all the calculations modulo $p$. We thus obtain relations with coefficients in the finite field $\mathbb{Z}/p\mathbb{Z}$. It is {\it a priori\/} possible that the rank of the relation matrix may become smaller (but never larger) after reduction modulo $p$. However, the `probability' of getting a wrong rank is very small (at most $1/p$, even smaller as the size of the matrix grows). Moreover, if the modulo $p$ rank is large enough, then the actual rank must also be large enough. In practice we also used different prime numbers to reduce the `risk'.

\subsection*{\it {\rm II. }Implementation remarks} --- During the calculation, even a single polynomial can be very big, so we use a linked list of monomials to store a polynomial. The order of the terms is not important, but when adding polynomials, the `combining like terms' procedure is crucial and must be done. To do this, we assign Hash values to the monomials, and use a Hash table to do the combination. This turns out to be the most efficient method.

Another way to speed up the calculation is to order the elements of ${\rm mon}_{(2g + 2, 2i)}$ properly, so that we can reach the expected dimensions by calculating the fewest relations possible. We do not know what the most efficient way of doing this is, but experimentally we found some good ordering of the relations, with no apparent theoretical reasons.

\smallskip
\noindent {\it The calculation was carried out on the Lisa Compute Cluster at the SARA Computing and Networking Services} (\url{http://www.sara.nl/}). {\it The program code is available from the authors upon request.}


\end{document}